\documentclass[12pt,oneside,english]{amsart}
\usepackage[T1]{fontenc}
\usepackage[latin1]{inputenc}
\makeatletter

\newcommand{\noun}[1]{\textsc{#1}}

 \theoremstyle{plain}
 \theoremstyle{plain}    
 \newtheorem{lem}{Lemma} 
 \theoremstyle{remark}
 \newtheorem{rem}{Remark}
 \theoremstyle{plain}    
 \newtheorem{thm}{Theorem} 
 \theoremstyle{definition}
  \newtheorem{example}{Example}
 \theoremstyle{plain}    
 \newtheorem{cor}{Corollary} 
 \theoremstyle{plain}    
 \newtheorem*{prop*}{Proposition} 

\usepackage{amsmath}
\usepackage{amsfonts}
\usepackage{amssymb}
\usepackage{babel}
\input{diagrams}
\input xy
\xyoption{all}

\def\CC{\mathbb{C}}

\def\QQ{\mathbb{Q}}
\def\QB{\bar{\QQ}}

\def\PP{\mathbb{P}}

\def\AA{\mathbb{A}}
\def\C{\mathcal{C}}

\def\D{\mathcal{D}}

\def\E{\mathcal{E}}
\def\F{\mathcal{F}}
\def\G{\mathcal{G}}
\def\H{\mathcal{H}}
\def\I{\mathcal{I}}
\def\J{\mathcal{J}}
\def\K{\mathcal{K}}
\def\L{\mathcal{L}}
\def\M{\mathcal{M}}

\def\N{\mathcal{N}}
\def\O{\mathcal{O}}

\def\T{\mathcal{T}}
\def\U{\mathcal{U}}

\def\V{\mathcal{V}}
\def\W{\mathcal{W}}

\def\Z{\mathcal{Z}}
\def\ts{\tilde{\ms}}
\def\tz{\tilde{\sZ}}
\def\btz{\bar{\tz}}
\def\bz{\bar{\sZ}}
\def\a{\alpha}
\def\b{\beta}
\def\d{\partial}
\def\g{\gamma}
\def\e{\epsilon}

\def\w{\omega}
\def\o{\Omega}
\def\l{\lambda}
\def\p{\text{p}}

\def\ms{\mathcal{S}}

\def\sV{\mathfrak{V}}
\def\sW{\mathfrak{W}}
\def\sX{\mathfrak{X}}

\def\sZ{\mathfrak{Z}}
\def\gs{\eta_{\ms}}

\def\rateq{\buildrel\text{rat}\over\equiv}

\def\homeq{\buildrel\text{hom}\over\equiv}
\def\nrateq{\mspace{7mu}/\mspace{-17mu}\rateq}
\def\nequiv{\mspace{7mu}/\mspace{-17mu}\equiv}

\def\bb{\bullet}
\def\mapa{\buildrel\theta\over\to}
\def\mapb{\buildrel\oplus \theta^{p,q}\over\to}

\def\mapf{\buildrel \cong \over \to}
\def\mapg{\buildrel \pi \over \to}
\def\mapgt{\buildrel \pi_{\T} \over \to}

\def\mapm{\buildrel \mathsf{X}_{\lambda} \over \to}

\def\mapphi{\buildrel \phi \over \to}
\def\mapbphi{\buildrel \overline{\phi} \over \to}
\def\mapdel{\buildrel c_{\D} \over \to}
\def\mapiota{\buildrel \iota \over \hookrightarrow}
\def\mapio{\buildrel \iota \over \to}
\def\mapnew{\buildrel '\iota^* \over \hookrightarrow}
\def\hm{Hom_{_{\text{MHS}}}}
\def\ext{Ext^1_{_{\text{MHS}}}}
\def\trdeg{\text{trdeg}}
\def\spec{\text{Spec}}
\def\ev{\text{ev}}

\def\pr{\text{pr}}
\def\un{\underline{\text{n}}}
\def\uf{\underline{\text{f}}}
\def\bv{\bar{\sV}}
\def\tv{\tilde{\sV}}
\def\btv{\bar{\tv}}

\usepackage{babel}
\makeatother
\begin{document}

\title{Exterior products of zero-cycles}

\author{Matt Kerr}

\begin{abstract}
We study the exterior product $CH_{0}(X)\otimes CH_{0}(Y)\to CH_{0}(X\times Y)$
on $0$-cycles modulo rational equivalence. The main tools used are
higher cycle- and $AJ$-classes developed in \cite{L1} and \cite{K2}.
The theorem of \cite{RS} (applied to $0$-cycles) appears as a special
case of our results.
\end{abstract}



\maketitle

\section{\textbf{Introduction}}

Since Jannsen's formal definition (in \cite{J4}) of a conjectural
\emph{Bloch-Beilinson filtration} on the Chow groups of smooth projective
varieties, a number of candidates have been put forth in the algebraic
geometry literature. Those of Murre \cite{M} and S. Saito \cite{sS}
are purely geometric, given in terms of the action of correspondences
on cycles. Raskind's approach \cite{Ra} is arithmetic, pulling back
a filtration on continuous \'etale cohomology (arising from the Hochschild-Serre
spectral sequence) along Jannsen's cycle-class map. On the other hand,
Griffiths-Green \cite{GG1}, Lewis \cite{L1} and M. Saito \cite{mS}
favor a Hodge-theoretic approach, using the Deligne-class map to pull
back a Leray filtration on cohomology to the Chow group; central here
is the idea of \emph{spreading out} a cycle.

Under reasonable conjectural assumptions, these filtrations not only
all yield BBF's --- they all coincide (e.g., see \cite{K1}). However,
they are still quite useful in the absence of these assumptions, for
instance in detecting cycles in the kernel of the Abel-Jacobi map.
Exterior products of homologically trivial cycles, yield such cycles;
and in this paper we turn our attention to the simplest case: products
$\Z_{1}\times\Z_{2}$ of degree-zero $0$-cycles, considered on the
product of the varieties on which they lie individually.

Let $Y_{1}$, $Y_{2}$ be smooth projective varieties of resp. dimensions
$d_{1}=1$, $d_{2}$; and let $\Z_{2}\in CH_{0}^{\text{{hom}}}(Y_{2})$
be such that $AJ(\Z_{2})$ is nontorsion in $J^{1}(Y_{2})$. By the
result of Rosenschon-Saito \cite{RS} one knows that if $Y_{1}$,
$Y_{2}$, $\Z_{2}$ are all defined over $K\subseteq\CC$ (say, finitely
generated $/\QB$) but $\Z_{1}\in CH_{0}^{\text{{hom}}}(Y_{1}/\CC)$
is not%
\footnote{i.e., up to rational equivalence class (a notion which makes perfect
sense as $CH_{0}(Y_{1}/\bar{K})\hookrightarrow CH_{0}(Y_{1}/\CC)$).%
} defined $/\bar{K}$, then $\Z_{1}\times\Z_{2}$ has infinite order
in $CH_{0}(Y_{1}\times Y_{2})$.%
\footnote{Below this will be written simply $\Z_{1}\times\Z_{2}$$\mspace{7mu}/\mspace{-21mu}\rateq0$,
because we take all cycle / cohomology groups $\otimes\QQ$ (i.e.,
modulo torsion). %
} The prototypical example is the $0$-cycle $(P_{1}-O_{1})\times(P_{2}-O_{2})$
on a product of curves $\C_{1}$,$\C_{2}$ defined $/\QB$, where
$P_{1}\in\C_{1}(\CC)$ is very general while $O_{1}\in\C_{1}(\QB)$,
$P_{2},O_{2}\in\C_{2}(\QB)$ and $AJ(P_{2}-O_{2})$ is nontorsion.

Our Theorem $1$ generalizes their result to $d_{1}>1$, replacing
the above condition on $\Z_{1}$ by (essentially) the requirement
that its $\bar{K}$-spread $\bar{\sZ_{1}}\in Z^{d_{1}}(Y_{1}\times\ms_{1})$
induce a nontrivial map of holomorphic forms $\o^{j}(Y_{1})\to\o^{j}(\ms_{1})$
for some $1\leq j\leq d_{1}$. Indeed, if $d_{1}=1$ and $\bar{\sZ_{1}}$
induces the zero map $\o^{1}(Y_{1})\to\o^{1}(\ms_{1})$, then one
can show that there exists a $\bar{K}$-specialization of $\Z_{1}$
having the same $AJ$ class, hence (as $d_{1}=1$ $\implies$ $Y_{1}$
is a curve) the same $\rateq$-class --- contradicting the \cite{RS}
condition on $\Z_{1}$. So for $d_{1}=1$ the conditions are equivalent.
The result leads to a generalization of the {}``prototypical example''
above to products of several curves; this was our original aim.

Returning to the Bloch-Beilinson picture, the filtration of Lewis
leads to the higher cycle- and Abel-Jacobi- classes $cl^{i}(\cdot)$
and $AJ^{i}(\cdot)$ of \cite{K2}; we describe how to compute these
below. The point is that if we interpret Theorem $1$ in terms of
these invariants, it says (modulo GHC) that $cl_{Y_{1}}^{j}(\Z_{1})\neq0$
and $AJ_{Y_{2}}^{(0)}(\Z_{2})[=Alb_{Y_{2}}(\Z_{2})]\neq0$ $\implies$
$AJ_{Y_{1}\times Y_{2}}^{j}(\Z_{1}\times\Z_{2})\neq0$, where the
latter $AJ^{j}$ is computed by a sort of cup product of $cl_{Y_{1}}^{j}$
and $AJ_{Y_{2}}$. This points the way to the much broader generalization
of Theorem $2$, where we start with nontrivial $higher$ invariants
$cl_{Y_{1}}^{j_{1}}(\Z_{1})$ and $AJ_{Y_{2}}^{j_{2}}(\Z_{2})$ for
both cycles and ask when the {}``cup product'' $AJ_{Y_{1}\times Y_{2}}^{j_{1}+j_{2}}(\Z_{1}\times\Z_{2})$
(corresponding to the exterior product of cycles) is nontrivial. The
answer involves a delicate quotient of the higher Abel-Jacobi class
$AJ_{Y_{2}}^{j_{2}}(\Z_{2})$ and careful consideration of the fields
of definition of $\Z_{1}$ and $\Z_{2}$. (In fact, there are two
different {}``quotients'' obtained by successive projections $AJ^{j}(\Z)\mapsto AJ^{j}(\Z)^{tr}\mapsto\overline{AJ^{j}(\Z)^{tr}}$;%
\footnote{note: since overlines are used to denote both quotients (as here)
and complex conjugation in this paper, footnotes usually alert the
reader in the latter case. (A bar also denotes completion / Zariski
closure but this is only for cycles and varieties.) %
} see eqn. $(4)$.)

Here is a more precise statement of our main results: let $Y_{1},\, Y_{2}$
be smooth projective (of any dimensions $d_{1},\, d_{2}$) and defined
$/\QB$, and $\L^{\bullet}$ denote Lewis's filtration on $CH_{0}$.
Denote by $\left\langle \Z\right\rangle \in CH_{0}$ the $\rateq$-class
of $\Z\in Z_{0}$.

\subsection*{Theorem $\mathbf{{1'}}$}

\emph{Given}\\
(a) \emph{a field $K\subseteq\CC$ finitely generated $/\QB$} (\emph{set
$j:=\trdeg(K/\QB)$})\emph{,}\\
(b) \emph{$\left\langle \Z_{1}\right\rangle \in\L^{j}CH_{0}(Y_{1}/K)$
with complete $\QB$-spread $\bar{\sZ_{1}}\in Z^{d_{1}}(Y_{1}\times\ms/\QB)$
inducing a nonzero map $\Omega^{j}(Y_{1})\to\Omega^{j}(\ms)$,}\\
(c) \emph{$\left\langle \Z_{2}\right\rangle \in CH_{0}^{hom}(Y_{2}/\QB)$
with nontorsion Albanese class in $Alb(Y_{2})$.}\\
\emph{Then $\Z:=\Z_{1}\times\Z_{2}\nrateq0$ in $\L^{j+1}CH_{0}(Y_{1}\times Y_{2}/K)$;
in particular, }\\
\emph{$AJ_{Y_{1}\times Y_{2}}^{j}(\Z)^{tr}\neq0$.}\\
\\
The two Corollaries provide various extensions --- to the case where
$j<\trdeg(K/\QB)$, or where a second finitely generated field $L$
takes the place of $\QB$.

Here is a simple consequence of this Theorem (from Example $2$ below).
If $\C_{i}/\QB$ ($i=1,\,\ldots,m+1$) are smooth projective curves
with $p_{i}\in\C_{i}(\CC)$ very general points, $o_{i}\in\C_{i}(\QB)$
($i=1,\,\ldots,m$), and $\W\in CH_{0}^{hom}(\C_{m+1}/\QB)$ has nontorsion
$AJ$-class in $J^{1}(\C_{m+1})$, then $(p_{1}-o_{1})\times\cdots\times(p_{m}-o_{m})\times\W\nrateq0$.
This was previously known only for $m=1$ (by \cite{RS}).

If Theorem $1'$ represents an application of new invariants to an
outstanding problem, the next Theorem can be seen as a statement about
the behavior of the invariants themselves under exterior product.

\subsection*{Theorem $\mathbf{{2'}}$}

\emph{Given}\\
(a) \emph{$K_{1},\, K_{2}\subseteq\CC$ f.g. $/\QB$ with compositum
$K$ satisfying $\trdeg(K_{1}/\QB)+\trdeg(K_{2}/\QB)=\trdeg(K/\QB)$,}\\
(b) \emph{$\left\langle \Z_{1}\right\rangle \in\L^{j_{1}}CH_{0}(Y_{1}/K_{1})$
with $cl_{Y_{1}}^{j_{1}}(\Z_{1})\neq0$,}\\
(c) \emph{$\left\langle \Z_{2}\right\rangle \in\L^{j_{2}}CH_{0}(Y_{2}/K_{2})$
with either}

(i) \emph{$cl_{Y_{2}}^{j_{2}}(\Z_{2})\neq0$ OR}

(ii) \emph{$\overline{AJ_{Y_{2}}^{j_{2}-1}(\Z_{2})^{tr}}\neq0$ and
$cl_{Y_{2}}^{j_{2}}(\Z_{2})=\cdots=cl_{Y_{2}}^{d_{2}}(\Z_{2})=0$.}\\
\emph{Assume the GHC.}\\
\emph{Then $\Z:=\Z_{1}\times\Z_{2}\nrateq0$ in $\L^{j_{1}+j_{2}}CH_{0}(Y_{1}\times Y_{2}/K)$.
In particular,}

\emph{when} (i) \emph{holds, $cl_{Y_{1}\times Y_{2}}^{j_{1}+j_{2}}(\Z)\neq0$;}

\emph{when} (ii) \emph{holds, $AJ_{Y_{1}\times Y_{2}}^{j_{1}+j_{2}-1}(\Z)^{tr}\neq0$.}\\
\\
The Proposition of $\S4$ states what can be proved without the GHC,
and has Theorem $1'$ as the special case corresponding to (ii) with
$j_{2}=1$.

In addition to these results there are several important lemmas which
will be valuable in further applications (e.g. \cite[sec. 7]{K1}).

This paper was written at the University of Chicago and MPI-Bonn;
we wish to thank both institutions for their hospitality.

\section{\textbf{Preliminaries}}

\subsection{Hodge structures}

A HS of weight m is a finite-dimensional $\QQ$-vector-space $\H_{(\QQ)}$
with a filtration $F^{\bb}$ on $\H_{\CC}:=\H\otimes_{\QQ}\CC$ such
that $F^{i}\H_{\CC}\oplus\overline{F^{m-i+1}\H_{\CC}}=\H_{\CC}=F^{0}\H_{\CC}$.
We denote $F^{i}\H_{\CC}\cap\overline{F^{m-i}\H_{\CC}}=:\H_{(\CC)}^{i,m-i}$
and note $\H_{\CC}=\oplus_{p+q=m}\H^{p,q}$. A $\QQ$-subspace $\G\subseteq\H$
is a subHS iff $\G_{\CC}=\oplus_{p+q=m}\G^{p,q}:=\oplus_{p+q=m}\left(\G_{\CC}\cap\H^{p,q}\right)$.
Intersections and sums of subHS are subHS. Moreover, the quotient
$\H/\G$ has a natural HS since $\H_{\CC}/\G_{\CC}\cong\oplus\left(\H^{p,q}/\G^{p,q}\right)$.
We write $F_{h}^{i}\H_{(\QQ)}$ for the largest subHS of $\H_{(\QQ)}$
contained%
\footnote{equality is only true in general for $m=2i$%
} in $F^{i}\H_{\CC}[\cap\H_{\QQ}]$, and $\QQ(-d)$ for the $1$-dimensional
weight $2d$ HS of pure type $(d,d)$.

For $\ms$ a smooth projective variety over a field $K\subseteq\CC$,
we write $H^{m}(\ms)$ for the HS $H_{\text{{sing}}}^{m}(\ms_{\CC}^{\text{{an}}},\QQ)$.
The fundamental class $[\Z]$ of an algebraic cycle%
\footnote{Note that for our purposes, this may be defined by integration and
Poincar\'e duality: if $\dim(\ms)=d$, then $\int_{\Z}(\cdot)\mapsto[\Z]$
under the identification $\{ H^{2d-2p}(\ms,\CC)\}^{\vee}\mapf H^{2p}(\ms,\CC)$.
It generates a subHS because it is in fact a rational $(p,p)$ class.%
} $\Z\in Z^{p}(\ms)$ gives a subHS $\QQ[\Z]\subseteq H^{2p}(\ms)$.
If $\ms=\ms_{1}\times\ms_{2}$, the direct summands under the K\"unneth
decomposition $H^{m}(\ms)=\oplus_{r+s=m}H^{r}(\ms_{1})\otimes H^{s}(\ms_{2})$
are subHS; and the K\"unneth components $[\Z]_{r}$ of $[\Z]$ give
subHS $\QQ[\Z]_{r}\subseteq H^{r}(\ms_{1})\otimes H^{2p-r}(\ms_{2})$.
(Of course tensor products of HS have natural HS.)

The Generalized Hodge Conjecture GHC($i,m,\ms$) predicts that $\linebreak$$F_{h}^{i}H^{m}(\ms)=N^{i}H^{m}(\ms)$,
where $N^{\bb}$ is the filtration by coniveau. The $F_{h}^{i}$ do
not behave so well under the K\"unneth decomposition: e.g., \[
F_{h}^{1}\left(H^{r}(\ms_{1})\otimes H^{s}(\ms_{2})\right)\supseteq F_{h}^{1}H^{r}(\ms_{1})\otimes H^{s}(\ms_{2})+H^{r}(\ms_{1})\otimes F_{h}^{1}H^{s}(\ms_{2})\]
 may be a proper inclusion. We will also need the following {}``skew''
subHS: \[
SF_{h}^{(i,j)}\left(H^{r}(\ms_{1})\otimes H^{s}(\ms_{2})\right)\,:=\,\text{{the\, largest\, subHS\, of\,}}H^{r}(\ms_{1})\otimes H^{s}(\ms_{2})\]
\[
\text{{contained\, in\,}}SF^{(i,j)}\left(H^{r}(\ms_{1},\CC)\otimes H^{s}(\ms_{2},\CC)\right)\cap\left[H^{r}(\ms_{1})\otimes H^{s}(\ms_{2})\right]\]
where\[
SF^{(i,j)}\left(H^{r}(\ms_{1},\CC)\otimes H^{s}(\ms_{2},\CC)\right):=\]
\[
F^{i}H^{r}(\ms_{1},\CC)\otimes F^{j}H^{s}(\ms_{2},\CC)+\overline{F^{i}H^{r}(\ms_{1},\CC)}\otimes\overline{F^{j}H^{s}(\ms_{2},\CC)}.\]
Note that $SF_{h}^{(1,\ell)}$ (for $\frac{s+1}{2}\geq\ell\geq0$,
$r>0$) contains \[
N^{1}H^{r}(\ms_{1})\otimes H^{s}(\ms_{2})+H^{r}(\ms_{1})\otimes F_{h}^{\ell}H^{s}(\ms_{2})\]
 (since $N^{1}\subseteq F_{h}^{1}\subseteq F_{\CC}^{1}\cap\overline{F_{\CC}^{1}}$,
$F^{\ell}H_{\CC}^{s}+\overline{F^{\ell}H_{\CC}^{s}}=H_{\CC}^{s}$,
etc.).

A morphism of HS $\H\mapa{}'\H$ is a $\QQ$-linear map which over
$\CC$ takes the form $\oplus\H^{p,q}\mapb\oplus{}'\H^{p,q}$ relative
to a pair of bases subordinate to the resp. Hodge decompositions.
Images and preimages of HS under such a morphism are HS.

We will use systematically the following notion of a {}``relative
dual pair'' of HS. This consists of: $\vspace{2mm}$\\
(a) $\K_{1}\subseteq\K_{0}$ of weights $2d-2n+1$, $\H_{1}\subseteq\H_{0}$
of weights $2n-1$, and$\vspace{2mm}$\\
(b) a perfect pairing $\H_{0}\times\K_{0}\to\QQ(-d)$ whose restriction
to $\H_{1}\times\K_{1}\to\QQ(-d)$ is also a perfect pairing.$\vspace{2mm}$\\
The inclusions (in (a)) induce (by the duality in (b)) projections
$\pr_{\H}:\,\H_{0}\twoheadrightarrow\H_{1}$, $\pr_{\K}:\,\K_{0}\twoheadrightarrow\K_{1}$,
and the compositions $\H_{1}\subseteq\H_{0}\twoheadrightarrow\H_{1}$,
$\K_{1}\subseteq\K_{0}\twoheadrightarrow\K_{1}$ yield the respective
identity maps.%
\footnote{One should view, for example, $\H_{1}\hookrightarrow\H_{0}$ as (a
choice of) extension of functionals from $\K_{1}$ to $\K_{0}$, and
$\H_{0}\twoheadrightarrow\H_{1}$ as restriction of functionals from
$\K_{0}$ to $\K_{1}$.%
} Equivalently, one could replace (b) with:$\vspace{2mm}$\\
(b$'$) identifications $\H_{i}\cong\K_{i}^{\vee}\otimes\QQ(-d)$
{[}which together with (a) induce $\pr_{\H}$, $\pr_{\K}${]} such
that the above compositions give the identity.$\vspace{2mm}$\\
In this paper the pairings always come tacitly from Poincar\'e duality.
Note that $\pr_{\H}$, $\pr_{\K}$ have kernels $\H_{1}'$, $\K_{1}'$
(resp.) which satisfy: $\H_{0}=\H_{1}\oplus\H_{1}'$, $\K_{0}=\K_{1}\oplus\K_{1}'$.
(This approach to complimentary HS's gives us more control than using
the semisimplicity coming from a polarization.) 

If $\H$ is of weight $2n-1$, define the Jacobian $J^{n}(\H):=\H_{\CC}/\left(F^{n}\H_{\CC}+\H_{\QQ}\right)$;
and if $\H=\K^{\vee}\otimes\QQ(-d)$, then $J^{n}(\H)=\left(F^{d-n+1}\K_{\CC}\right)^{\vee}/\K_{\QQ}^{\vee}$.
When $\H_{1}\subseteq\H_{0}$ is a subHS, $J^{n}(\H_{1})\hookrightarrow J^{n}(\H_{0})$
and $J^{n}(\H_{0})\twoheadrightarrow J^{n}(\H_{0}/\H_{1})\cong J^{n}(\H_{0})/J^{n}(\H_{1})$.
In the above {}``dual pair'' situation, $J(\H_{0})=J(\H_{1})\oplus J({}'\H_{1})$;
we emphasize that since extension followed by restriction of functionals
$\H_{1}\subseteq\H_{0}\twoheadrightarrow\H_{1}$ is the identity,
so is the induced composition on Jacobians. More generally a morphism
$\theta$ induces a map of Jacobians. We write $J^{p}(H^{2p-1}(\ms))=:J^{p}(\ms)$.

To construct elements in Jacobians: let $\dim(\ms)=d$, $\Z\in Z_{\text{{hom}}}^{p}(\ms)$,
$\d^{-1}\Z\,=$ any choice of topological $(2d-2p+1)$-chain bounding
on $\Z$, and let $\K\subseteq H^{2d-2p+1}(\ms)$, $\H\subseteq H^{2p-1}(\ms)$
be a relative dual pair. Then $\int_{\d^{-1}\Z}(\cdot)$ is a well-defined
functional on $F^{d-p+1}H^{2d-2p+1}(\ms,\CC)$. (Represent the latter
by $C^{\infty}$ forms $\w\in F^{d-p+1}\o_{\ms^{\infty}}^{2d-2p+1}(\ms)$;
note that if $\w-\w'=d\a$, $\a$ may be chosen $\in F^{d-p+1}$,
and use Stokes's theorem.) Hence one has an element of $\left\{ F^{d-p+1}\K_{\CC}\right\} ^{\vee}\twoheadrightarrow J^{p}(\H)$;
we write $\left\langle \int_{\d^{-1}\Z}(\cdot)\right\rangle \in J^{p}(\H)$.

\subsection{The fundamental lemma}

Let $\H_{1}\subseteq\H_{0}$, $\K_{1}\subseteq\K_{0}$ be a relative
dual pair of HS, with $\G_{0}\subseteq\H_{0}$ a subHS closed under
$\H_{0}\twoheadrightarrow\H_{1}\subseteq\H_{0}$. Furthermore, let
$\H_{\sV}\subseteq\H_{1}$, $\K_{\sV}\subseteq\K_{1}$ be another
such pair, with $\H_{2}$ and $\K_{2}$ the resp. projection-kernels:
$\H_{1}=\H_{\sV}\oplus\H_{2}$, etc. One has the following diagram:

\xymatrix{ & { \frac{ \{ F^{d-n+1}\K_0^{\CC} \}^{\vee}}{ \{ \K_0^{\QQ} \}^{\vee}} } \ar @{<->} [r]^{\mspace{30mu}\cong} & {J(\H_0)} \ar @{>>} [d]^{\pr_1} \ar @{>>} [r]^{\b_0} & {J \left( \frac{\H_0}{\G_0} \right) } \ar @{>>} [d]^{\overline{\pr_1}} \\ {\frac{ \{ F^{d-n+1} \K_{\sV}^{\CC} \}^{\vee}}{ \{ \K_{\sV}^{\QQ} \}^{\vee} }} \ar @{<->} [r]^{\cong} & {J(\H_{\sV})} \ar @{^{(}->} [r]^{\iota_{\sV}} & {J(\H_1)} \ar [d]^{(\pr_{\sV},\pr_2)}_{\cong} \ar @{>>} [r]^{\b_1} & {J \left( \frac{\H_1}{\H_1 \cap \G_0} \right) } \ar @{=} [d] \\ & & {J(\H_{\sV})\oplus J(\H_2)} \ar @{>>} [ul] & {\frac{J(\H_1)}{J(\H_1 \cap \G_0 ) } } }

${}$\\
in which the square commutes and $\pr_{\sV}\circ\iota_{\sV}$ is the
identity.

Now take $\Xi\in J(\H_{0})$ with lifting $\tilde{\Xi}\in\{ F^{d-n+1}\K_{0}^{\CC}\}^{\vee}$,
$\xi\in J(\H_{\sV})$ with lifting $\tilde{\xi}\in\{ F^{d-n+1}\K_{\sV}^{\CC}\}^{\vee}$,
such that $\tilde{\Xi}$ is trivial on $F^{d-n+1}\K_{2}^{\CC}$ $[=\K_{2}^{\CC}\cap F^{d-n+1}\K_{0}^{\CC}]$
and equivalent to $\tilde{\xi}$ on $F^{d-n+1}\K_{\sV}^{\CC}$. Then
$(\pr_{\sV}\circ\pr_{1})(\Xi)=\xi$, $(\pr_{2}\circ\pr_{1})(\Xi)=0$,
hence $\pr_{1}(\Xi)=\iota_{\sV}(\xi)$.

\subsection*{Case 1:}

Suppose $\xi\neq0$ and $\H_{\sV}\cap(\H_{1}\cap\G_{0})=\{0\}$. Then
$J(\H_{\sV})\oplus J(\H_{1}\cap\G_{0})\hookrightarrow J(\H_{1})$
and hence $\b_{1}$ cannot kill $\iota_{\sV}(\xi)$; we conclude that
$\beta_{0}(\Xi)\neq0$.

\subsection*{Case 2:}

Let $\G_{1}\subseteq\H_{1}$ and suppose that $\pi_{\sV}(\xi)\neq0$
in the diagram

\xymatrix{ {J(\H_\sV)} \ar @{>>} [d]^{\pi_{\sV}} \ar @{^{(}->} [r]^{\iota_\sV } & {J(\H_1)} \ar @{>>} [d]^{\pi} \ar @{>>} [r]^{\b_1} & {J \left( \H_1 / (\H_1 \cap \G_0 ) \right) } \ar @{>>} [d]^{\overline{\pi}} \\ {J \left( \H_{\sV} /(\G_1 \cap \H_{\sV}) \right) } \ar @{^{(}->} [r]^{\mspace{40mu} \overline{\iota_{\sV}}} & {J \left( \H_1 / \G_1 \right) } \ar @{>>} [r]^{\overline{\b_1} \mspace{70mu}} & {J \left( \H_1 / \{ (\H_1 \cap \G_0 ) + \G_1 \} \right) } \ar @{=} [d] \\ & & {\frac{J \left( \H_1 / \G_1 \right)}{J \left( ( \H_1 \cap \G_0) / (\G_1 \cap \G_0 ) \right) } } }${}$\\
(in which squares commute). Furthermore assume only that $(\H_{1}\cap\G_{0})\cap\H_{\sV}\subseteq\G_{1}$.
This implies that $\frac{\H_{1}\cap\G_{0}}{\G_{1}\cap\G_{0}}\cap\frac{\H_{\sV}}{\G_{1}\cap\H_{\sV}}=\{0\}$
in $\frac{\H_{1}}{\G_{1}}$, hence $J\left(\frac{\H_{1}\cap\G_{0}}{\G_{1}\cap\G_{0}}\right)\oplus J\left(\frac{\H_{\sV}}{\G_{1}\cap\H_{\sV}}\right)\hookrightarrow J\left(\frac{\H_{1}}{\G_{1}}\right)$.
Therefore $\overline{\b_{1}}$ cannot kill $\overline{\iota_{\sV}}(\pi_{\sV}(\xi))$,
and so (by both diagrams) $\overline{\pi}(\overline{\pr_{1}}(\b_{0}(\Xi)))=\overline{\b_{1}}(\overline{\iota_{\sV}}(\pi_{\sV}(\xi)))\neq0$
$\implies$ $\b_{0}(\Xi)\neq0$ once again.

\subsection{Points and spreads}

Given $\ms/\QB$ smooth projective, choose any affine Zariski open
subset $S/\QB$. The embedding $\QB[S]\hookrightarrow\QB(S)\cong\QB(\ms)$
then produces a generic point $p_{g}$ on $\ms$ via the composition$\linebreak$$\spec\,\QB(\ms)\to\spec\,\QB[S]\,\cong S\hookrightarrow\ms$.
For purposes of taking cohomology, we use the approximation $\gs=\underleftarrow{\lim}\,\U$
(over $\U\subseteq\ms$ affine Zariski open subsets defined $/\QB$)
to $p_{g}$; more precisely, $H^{i}(\gs):=\underrightarrow{\lim}\, H^{i}(\U_{\CC}^{\text{{an}}},\QQ)$
while $CH^{p}(X\times_{\QB}\gs):=\underrightarrow{\lim}\, CH^{p}(X\times_{\QB}\U)\cong CH^{p}(X\times_{\QB}p_{g})$.
Note that $\underline{H^{i}}(\gs):=W_{i}H^{i}(\gs)=\text{{im}}\left\{ H^{i}(\ms)\to H^{i}(\gs)\right\} \cong H^{i}(\ms)/N^{1}H^{i}(\ms)$
is a HS.%
\footnote{In general for quasiprojective $Y$, we write $\underline{H^{i}}(Y):=\text{{im}}\{ H^{i}(\overline{Y})\to H^{i}(Y)\}$;
this is independent of the choice of smooth compactification $Y$.%
}

Given any embedding $\text{{ev}}:\,\QB(\ms)\hookrightarrow\CC$ which
restricts to the indentity on $\QB$, define a geometric point $p\in\ms(\CC)$
of maximal transcendence degree ($=\dim(\ms)$) over $\QB$ by $p:=p_{g}\times_{\text{{ev}}}\spec\,\CC$.
(See \cite{K1}. In fact $p$ is defined $/K$, for $K:=\text{{ev}}(\QB(\ms))$.)
Since $p$ lies in the complement of the countably many divisors $D\subseteq\ms$
defined $/\QB$, we may think of it as a geometric (closed, zero-dimensional)
point of $\gs^{\CC}$; such a point is called \emph{very general}.

To define the $\QB$-spread of an algebraic cycle we need part $(a)$
of the following:

\begin{lem}
(a) Let $K\subseteq\CC$ be a finitely generated extension of $\QB$.
Then $\exists$ $\ms/\QB$ smooth projective and a very general point
$p\in\ms(\CC)$ such that $\text{{ev}}_{p}:\,\QB(\ms)\mapf K$.

(b) Given $\QB(\ms_{1})\cong K_{1}$, $\QB(\ms_{2})\cong K_{2}$ two
such, $\QB(\ms_{1}\times\ms_{2})\cong\QB(K_{1},K_{2})$ $\Longleftrightarrow$
$\trdeg(\QB(K_{1},K_{2})/\QB)=\trdeg(K_{1}/\QB)+\trdeg(K_{2}/\QB)$.
\emph{{[}}Here $\QB(K_{1},K_{2})$ denotes the compositum of $K_{1}$
and $K_{2}$.\emph{{]}}
\end{lem}
\begin{proof}
$(a)$ $K=\QB(\pi_{1},\ldots,\pi_{t};\,\a_{1},\ldots,\a_{s})$ for
$\{\pi_{i}\}$ a transcendence basis $\implies$ we can map $\QB[x_{1},\ldots,x_{t};\, x_{t+1},\ldots,x_{t+s}]\mapphi K$
via $\{ x_{i}\mapsto\pi_{i},\, x_{t+j}\mapsto\a_{j}\}$. Set $R=\text{{im}}(\phi)$,
$I=\ker(\phi)$, $S=\text{{Var}}(I)\subseteq\AA_{\QB}^{t+s}$. Since
$R$ is a division ring, $I$ is prime and $S$ irreducible, $\dim S=t=\trdeg(K/\QB)$.
So $\overline{\phi}:\,\QB[S]=\frac{\QB[x_{1},\ldots,x_{t+s}]}{I}\mapf R$
induces an isomorphism of fraction fields%
\footnote{$\widetilde{S}$ denotes a desingularization of $S$.%
} $\QB(S)=\frac{\QB(x_{1},\ldots,x_{t})[x_{t+1},\ldots,x_{t+s}]}{I}\mapf K$;
this is evaluation at $p=(\pi_{1},\ldots,\pi_{t};\,\a_{1},\ldots,\a_{s})\in S(\CC)$.
Finally, take $\widetilde{S}$ to be a desingularization of $S$,
and $\ms$ to be a good compactification of $\widetilde{S}$; one
has $\QB(\ms)\cong\QB(\widetilde{S})\cong\QB(S)$.

$(b)$ As in $(a)$, we have $p=p_{1}\times p_{2}=(\{\underline{\pi}\},\{\underline{\a}\};\{\underline{\sigma}\},\{\underline{\b}\})\in S_{1}\times S_{2}\subseteq\AA_{\{ x\}}^{t_{1}+s_{1}}\times\AA_{\{ y\}}^{t_{2}+s_{2}}$
for $\{\underline{\pi}\},\{\underline{\sigma}\}$ transcendence bases
for $K_{1},K_{2}\,/\QB$. Evaluation at $p$ gives a map\SMALL\[
\QB[S_{1}\times S_{2}]=\QB[S_{1}]\otimes_{\QB}\QB[S_{2}]\cong\frac{\QB[x_{1},\ldots,x_{t_{1}+s_{1}};y_{1},\ldots,y_{t_{2}+s_{2}}]}{I=(I_{1}^{\{ x\}},I_{2}^{\{ y\}})}\mapbphi\QB(K_{1},K_{2})\subseteq\CC.\]
\normalsize If $\QB(\ms_{1}\times\ms_{2})[\cong\QB(S_{1}\times S_{2})]$
is not $\cong$ to the fraction field of $\text{{im}}(\overline{\phi})$,
then $\overline{\phi}$ kills some $f\nequiv0$ (mod $I$). Since
$I$ is prime, Nullstellensatz $\implies$ $\I(S_{1}\times S_{2})=I$,
hence $f$ does not vanish on $S_{1}\times S_{2}$ and $f=0$ cuts
out a subvariety $D/\QB$ of codim.$\geq1$ in which $p$ must sit.
Since the relative dimension of $\pr:\, S_{1}\times S_{2}\to\AA_{\{ x\}}^{t_{1}}\times\AA_{\{ y\}}^{t_{2}}$
is $0$, $\pr(p)=(\{\underline{\pi}\},\{\underline{\sigma}\})$ sits
in a $\QB$-subvariety of $\AA^{t_{1}+t_{2}}$ of codim.$\geq1$;
therefore $\pi_{1},\ldots,\pi_{t_{1}};\,\sigma_{1},\ldots,\sigma_{t_{2}}$
are not algebraically independent. But $\{\a_{i},\b_{j}\}$ are algebraic
over $\QB(\{\pi\},\{\sigma\})$, hence $\QB(K_{1},K_{2})$ has $\trdeg\leq t_{1}+t_{2}-1$.
\end{proof}
Now let $X$ be defined over $\QB$, $K$ be f.g. $/\QB$, $X_{K}=X\otimes_{\QB}\spec\, K$
and $\Z\in Z^{p}(X_{K})$. (All cycle and Chow groups will be taken
$\otimes\QQ$; we also write $Z^{*}(X/K)$ for $Z^{*}(X_{K})$.) By
Lemma $1(a)$ one has $\ev_{p}:\,\QB(\ms)\mapf K$, and we define
$\Z_{g}:=\Z\times_{\ev^{-1}}\spec\,\QB(\ms)\in Z^{p}(X_{\QB(\ms)})$.
Clearing denominators from the equations cutting out the components
of $\Z_{g}$ yields a cycle $\bar{\sZ}\in Z^{p}(X\times\ms/\QB)$
whose complexification restricts to $\Z$ along $X\times\{ p\}\hookrightarrow X\times\ms$.
(Such a {}``complete spread'' is not well-defined modulo $\rateq$.)
One also has the obvious restriction $\sZ\in Z^{p}(X\times\gs/\QB)$
which is called the $\QB$-spread of $\Z$. We can write this as a
map

\begin{equation} CH^p(X_K) \mapf CH^p(X\times \gs /\QB ) . \end{equation}${}$

\subsection{Higher cycle- and Abel-Jacobi classes}

To take the Deligne class of the r.h.s. of $(1)$ we use the well-defined
{}``composition'' $\psi$:\small\[
CH^{p}(X\times\gs/\QB)\twoheadleftarrow CH^{p}(X\times\ms/\QB)\mapdel H_{\D}^{2p}(X\times\ms,\QQ(p))\twoheadrightarrow\underline{H_{\D}^{2p}}(X\times\gs,\QQ(p))\]
\normalsize (where $\underline{H_{\D}^{2p}}$ denotes the image of
Deligne cohomology of $(X\times\ms)_{\CC}^{\text{{an}}}$ in absolute
Hodge cohomology of $(X\times\gs)_{\CC}^{\text{{an}}}$, see \cite{L1}).
Write $\Psi^{K/\QB}:=\psi\circ(1)$.

Lewis \cite{L1} constructs a Leray filtration $\L^{\bb}$ on $\underline{H_{\D}^{2p}}$
where (in the notation of \cite{K2})

\small \begin{equation} 0 \to Gr^{i-1}_{\L} \underline{J^p} (X\times \gs ) \to Gr^i_{\L} \underline{H^{2p}_{\D}} (X\times \gs, \QQ(p) ) \to Gr^i_{\L} \underline{Hg^p} (X\times \gs ) \to 0 . \end{equation} \normalsize Taking
$\psi$-preimages gives $\L$ on both groups of $(1)$; and if $\sZ\in\L^{i}$($\Longleftrightarrow\Z\in\L^{i}$)
then its invariants in $(2)$ are written $[c_{\D}(\sZ)]_{i}$ (or
$\Psi_{i}^{K/\QB}(\Z)$), $[\sZ]_{i}$, and (if $[\sZ]_{i}=0$) $[AJ(\sZ)]_{i-1}$.
One has

\begin{equation} \begin{matrix} Gr^i_{\L} \underline{Hg^p} (X \times \gs) \cong \hm \left(  \QQ(-p) , \underline{H^i}(\gs) \otimes H^{2p-i}(X) \right) \\ {} \\ \hookrightarrow \underline{H^i}(\gs) \otimes H^{2p-i}(X) , \end{matrix} \end{equation}${}$

\begin{equation} \begin{matrix} Gr^{i-1}_\L \underline{J^p}(X\times \gs) \cong \frac{\ext \left( \QQ(-p), \underline{H^{i-1}}(\gs) \otimes H^{2p-i}(X) \right) }{\text{im}\, \{ \hm \left( \QQ(-p), Gr^W_i H^{i-1}(\gs) \otimes H^{2p-i}(X) \right) \} } \\ {} \\ \twoheadrightarrow J^p \left( \underline{H^{i-1}}(\gs) \otimes \frac{H^{2p-i}(X)}{F^{p-i+1}_h H^{2p-i}(X)} \right) \twoheadrightarrow J^p \left( \frac{H^{i-1}(\ms)\otimes H^{2p-i}(X)}{SF^{(1,p-i+1)}_h \{\text{num}\}} \right) , \end{matrix} \end{equation}${}$\\
where the first projection is worked out in \cite[sec. 12]{K2}.%
\footnote{The second projection (where {}``num'' just means {}``numerator'')
is valid by the remarks on $SF$ in $\S2.2$ for $i\geq2$, and trivially
for $i=1$.%
} The projected images of $[AJ(\sZ)]_{i-1}$ are written $[AJ(\sZ)]_{i-1}^{tr}$
and $\overline{[AJ(\sZ)]_{i-1}^{tr}}$.%
\footnote{where the bar denotes image under projection, not complex conjugation.%
} (Note that the $Ext\cong J^{p}(\underline{H^{i-1}}(\gs)\otimes H^{2p-i}(X))$.)

To compute $[\sZ]_{i}$ one takes the image of the K\"unneth component
$[\bar{\sZ}]_{i}\in H^{i}(\ms)\otimes H^{2p-i}(X)$ in the r.h. term
of $(3)$. If $[\bar{\sZ}]=0$ then one may compute $[AJ(\sZ)]_{i-1}$
and its images by projecting $\left\langle \int_{\d^{-1}\bar{\sZ}}(\cdot)\right\rangle \in J^{p}(H^{i-1}(\ms)\otimes H^{2p-i}(X))$
to the appropriate term in $(4)$.

Taking a limit over $K\subseteq\CC$ f.g. $/\QB$ (and using $CH^{p}(X_{K})\hookrightarrow CH^{p}(X_{\CC})$
and the corresponding limit%
\footnote{i.e. the only maps $Gr_{\L}^{i}\underline{H_{\D}^{2p}}(X\times\gs,\QQ(p))\to Gr_{\L}^{i}\underline{H_{\D}^{2p}}(X\times\eta_{\ms'},\QQ(p))$
allowed in the definition of the limit are those corresponding to
$dominant$ morphisms $\ms'\to\ms$; see \cite[sec. 3]{K2}.%
} of $\ms$'s) we have filtrations and maps\[
\Psi_{i}^{(\CC/\QB)}:\,\L^{i}CH^{p}(X_{\CC})\to\underrightarrow{\lim}\, Gr_{\L}^{i}\underline{H_{\D}^{2p}}(X\times\gs,\,\QQ(p))\,;\]
write invariants $cl_{X}^{i}(\Z)$ and (if this $=0$) $AJ_{X}^{i-1}(\Z)$
(which are essentially $[\sZ]_{i}$ and $[AJ(\sZ)]_{i-1}$ but without
referring to $K$ or $\sZ$).

We will use the following two lemmas in the proofs of Corollary $1$
and Theorem $2$. In the first one (writing $d=\dim X$) we take $p=d$,
which corresponds to the case where $\Z$ is a $0$-cycle.

\begin{lem}
Let $\sZ\in\L^{i}CH^{d}(X\times\gs/\QB)$. If $[\sZ]_{i}\neq0$ and
GHC$(1,i,\ms)$ holds then the induced map $\bar{\sZ}^{*}:\,\o^{i}(X)\to\o^{i}(\ms)$
is nontrivial.
\end{lem}
\begin{proof}
Noting that $N^{1}H^{i}(\ms)\subseteq F_{h}^{1}H^{i}(\ms)$ $\implies$
$\underline{H^{i,0}}(\gs,\CC)=H^{i,0}(\ms,\CC)$, we must show the
composition $\phi:$\[
\hm\left(\QQ(-p),\underline{H^{i}}(\gs)\otimes H^{2d-i}(X)\right)\hookrightarrow\underline{H^{i}}(\gs,\CC)\otimes H^{2d-i}(X,\CC)\]
\[
\twoheadrightarrow H^{i,0}(\ms,\CC)\otimes H^{d-i,d}(X,\CC)\]
is injective. An element of the $\hm$ is a morphism $\theta:\, H^{i}(X)\to\underline{H^{i}}(\gs)$
of HS; if $\theta\in\ker(\phi)$ then \[
\text{{im}}(\theta)\subseteq\ker\left\{ \underline{H^{i}}(\gs)\hookrightarrow\underline{H^{i}}(\gs,\CC)\twoheadrightarrow H^{i,0}(\ms,\CC)[\oplus H^{0,i}(\ms,\CC)]\right\} .\]
Since moreover $\text{{im}}(\theta)$ is a HS, we have $\text{{im}}(\theta)\subseteq F_{h}^{1}\underline{H^{i}}(\gs)$;
and GHC $\implies$ $F_{h}^{1}\underline{H^{i}}(\gs)\subseteq N^{1}\underline{H^{i}}(\gs)=0$.
\end{proof}
${}$

\begin{lem}
(a) For $\ts/\QB$ smooth projective of dimension $i+c$, $\exists$
$c$-fold hyperplane section $\ms/\QB$ s.t. restriction along $\ms\mapiota\ts$
induces a well-defined injection $\underline{H^{i}}(\eta_{\ts})\hookrightarrow\underline{H^{i}}(\gs)$.

(b) Given $\btz\in Z^{p}(X\times\ts/\QB)$ with restriction $\bz$
to $X\times\ms$ (and write $\tz,\,\sZ$ for their resp. restrictions
to $X\times\eta_{\ts},\, X\times\gs$). If $c_{\D}(\tz),\, c_{\D}(\sZ)$
belong to $\L^{i}$ of the resp. $\underline{H_{\D}^{2p}}$'s, then
$\iota$ induces a (well-defined) injective map of invariants $[\tz]_{i}\mapsto[\sZ]_{i}$.
If they belong to resp. $\L^{i+1}$'s and $[\btz]=0$ $[\implies[\bz]=0]$,
then $\iota$ induces injections $[AJ(\tz)]_{i}^{tr}\mapsto[AJ(\sZ)]_{i}^{tr}$,
$\overline{[AJ(\tz)]_{i}^{tr}}\mapsto\overline{[AJ(\sZ)]_{i}^{tr}}$.

(c) Assume HC. Then if $c_{\D}(\tz)\in\L^{j}\underline{H_{\D}^{2p}}(X\times\eta_{\ts},\QQ(p))$,
one can choose $\btz$ so that (in (b)) $c_{\D}(\sZ)\in\L^{j}\underline{H_{\D}^{2p}}(X\times\gs,\QQ(p))$.
\end{lem}
\begin{proof}
(a) Arguing for $c=1$, let $\ms\subseteq\tilde{\ms}$ be a smooth
$\QB$-hyperplane section and choose resp. codim.-$1$ $\QB$-subvarieties
$D\subseteq\ms$, $\tilde{D}\subseteq\tilde{\ms}$ as follows: $D$
sufficiently {}``large'' that $\underline{H^{i}}(\ms\setminus D)\mapf\underline{H^{i}}(\gs)$;
and $\tilde{D}$ properly intersecting $\ms$ with $\tilde{D}\cap\ms\supseteq D$.
By \cite[Thm. 6.1.1]{AS} (a version of affine weak Lefschetz), $H^{i}(\tilde{\ms}\setminus\tilde{D})\hookrightarrow H^{i}(\ms\setminus\tilde{D}\cap\ms$)
and so $\underline{H^{i}}(\tilde{\ms}\setminus\tilde{D})\mapnew\underline{H^{i}}(\ms\setminus\tilde{D}\cap\ms)$.
Moreover, by \cite[Thm. 1.1(3)]{AK}, $\iota^{*}:\, H^{i}(\tilde{\ms})\hookrightarrow H^{i}(\ms)$
respects the coniveau filtration. Hence we get a commutative diagram

\begin{equation*} \xymatrix{{H^i(\tilde{\ms})} \ar @{>>} [r] \ar @{^{(}->} [d]_{\iota^*} & {\frac{H^i(\tilde{\ms})}{\text{im}{H^i_{\tilde{D}}(\tilde{\ms})}}} \ar @{>>} [r] \ar @{^{(}->} [d]_{'\iota^*} & {\frac{H^i(\tilde{\ms})}{N^1H^i(\tilde{\ms})}} \ar @{<->} [r]^{\cong} \ar [d]_{''\iota^*} & {\underline{H^i}(\eta_{\tilde{\ms}})} \\ H^i(\ms) \ar @{>>} [r] & {\frac{H^i(\ms)}{\text{im}{H^i_{\tilde{D}\cap\ms}(\ms)}}} \ar [r]^{\cong} & {\frac{H^i(\ms)}{N^1H^i(\ms)}} \ar @{<->} [r]^{\cong} & \underline{H^i}(\gs) } \end{equation*}${}$\\
from which injectivity of $''\iota^{*}$ is obvious. (Iterating this
procedure proves it for $c>1$.)

(b) The point is to plug the injection of (a) into $\linebreak$$\hm\left(\QQ(-p),\,(\text{{---}})\otimes H^{2p-i}(X)\right)$,
$J^{p}\left((\text{{---}})\otimes\frac{H^{2p-i-1}(X)}{F_{h}^{p-i}H^{2p-i-1}(X)}\right)$,
and $\linebreak$$J^{p}\left(\frac{(\text{{---}})\otimes H^{2p-i-1}(X)}{SF_{h}^{(1,p-i)}\{\text{{num}}\}}\right)$.
This automatically yields injections except in the last case, where
we need $\frac{\underline{H^{i}}(\eta_{\ts})\otimes H^{2p-i-1}(X)}{SF_{h}^{(1,p-i)}}\hookrightarrow\frac{\underline{H^{i}}(\gs)\otimes H^{2p-i-1}(X)}{SF_{h}^{(1,p-i)}}$
(in order that the $J^{p}$'s inject). It suffices to show \[
SF_{h,\,\ms}^{(1,p-i)}\cap\text{{im}}\left\{ \underline{H^{i}}(\eta_{\tilde{\ms}})\otimes H^{2p-i-1}(X)\right\} \,=\, SF_{h,\,\ts}^{(1,p-i)}.\]
 This follows by describing $SF_{h}^{(1,p-i)}$ as the largest subHS
contained in $\ker\left\{ \underline{H^{i}}(\eta)\otimes H^{2p-i-1}(X)\to\underline{H^{i,0}}(\eta,\CC)\otimes H^{p-i-1,\, p}(X,\CC)\right\} $
and noting that $\underline{H^{i,0}}(\eta_{\ts},\CC)\hookrightarrow\underline{H^{i,0}}(\gs,\CC)$.

(c) By \cite[sec. 1.6]{mS} we can arrange that $c_{\D}(\btz)\in\L_{(X\times\ts)/\ts}^{i}H_{\D}^{2p}(X\times\ts,\QQ(p))$,
by modifying an initial choice of $\btz$ along $X\times D$ where
$D\subseteq\ts$ is a divisor defined $/\QB$. (This uses the HC.)
The conclusion then follows by functoriality of Leray.
\end{proof}
\begin{rem}
(i) We emphasize that the passage from $\tz\mapsto\btz\mapsto\bz\mapsto\sZ$
is not a well-defined map of cycles (only $\btz\mapsto\bz$ is). However,
(b) says that at least certain $invariants$ of $\sZ$ will only depend
on those of $\tz$, and not on the choice of lifting $\tz\mapsto\btz$.

(ii) Given $\tilde{\Z}$ with (complete) spread $\btz$, $\bz$ should
be viewed as a spread of a {}``specialization'' $\Z$ of $\tilde{\Z}$
(defined over a field of lesser transcendence degree $/\QB$). Namely,
if $p\in\ms(\CC)$ is very general (hence somewhat $less$ than very
general in $\ts(\CC)$), take $\Z$ to be the restriction of $\bz$
along $X\times\{ p\}\hookrightarrow X\times\ms$. {[}Note that this
is different (less delicate) than the kind of specialization considered
in \cite{GGP}, \cite{mS}.{]}
\end{rem}

\subsection{Change of spread field}

We now make a slight extension to the case where $X$ is not defined
$/\QB$, to be used in Corollary $2$ (and, to a lesser extent, Example
$3$).

Suppose we have $L\subseteq K\,[\subseteq\CC]$ both f.g. $/\QB$,
with $\trdeg(K/L)=:t\geq1$. Then $\exists$ $\ms/\QB$ with $s\in\ms(\CC)$
such that $\ev_{s}:\,\QB(\ms)\mapf K$, and $\M/\QB$ with a morphism%
\footnote{If $\ms$, $\M$ come from Lemma $1(a)$, $\rho$ is \emph{a priori}
a rational map, restricting to a morphism only on $U\subseteq\ms$
Zariski open. Take the closure in $\ms\times\M$ of $graph(\rho|_{U})\subseteq U\times\M\subseteq\ms\times\M$
and let $\ms'$ be a desingularization of the result. Then one has
obvious morphisms $\ms'\twoheadrightarrow\ms$ and $\ms'\to\M$, and
the first is a birational equivalence $(\QB(\ms')\cong\QB(\ms)$);
so just take {}``$\ms$'' in the above to be $\ms'$.%
} $\rho:\,\ms\to\M$ such that $\ev_{\rho(s)}:\,\QB(\M)\mapf L$. Write
$\mu=\rho(s)\in\M(\CC)$ and $\T:=\rho^{-1}(\mu)\mapiota\ms_{L}$,
and note that $L(\T)\cong K$ (again via $\ev_{s}$). Let $X$ be
defined $/L$ and $\Z\in Z^{p}(X_{K})$; then one has complete $\QB$-spreads
$\bz\in Z^{p}(\bar{\sX}/\QB)$ and $\bar{\sX}\mapg\ms$, formally
restricting to $\sX:=\underleftarrow{\lim}\,\pi^{-1}(\U)$ {[}over
$\U_{/\QB}\subseteq\ms$ affine Zar. op.{]} and $\sZ\in Z^{p}(\sX)$.
Moreover, one has the \emph{partial} ($L$-)spreads\[
\bar{\sX}_{L}\times_{\ms_{L}}\T=X_{L}\times\T\mapiota\bar{\sX}_{L}\text{{\,\, and\,\,}}\bar{\sZ_{\T}}:=\iota^{*}(\bz_{L})\in Z^{p}(X_{L}\times\T)\,;\]
we write $\eta_{\T}=\underleftarrow{\lim}\,\V$ {[}$\V_{/L}\subseteq\T$
affine Zar. op.{]} and note that this is formally the restriction
to $\T$ of $(\gs)_{L}$.

By functoriality of the Deligne class we get a commuting diagram

\xymatrix{ {CH^p(X/K)} \ar [r]^{\cong} \ar [dr]_{\cong} & {CH^p(\sX /\QB) \ar [r]^{c_{\D}} } \ar [d] & {\underline{H^{2p}_{\D}} (\sX^{\text{an}}_{\CC} ,\QQ(p) ) } \ar [d]^{\iota^*} \\ & {CH^p(X_L \times \eta_{\T} /L )} \ar [r]^{c_{\D} \mspace{25mu}} & {\underline{H^{2p}_{\D}} \left( (X \times \eta_{\T} \right) ^{\text{an}}_{\CC} ,\QQ(p))} } ${}$\\
where the two $\cong$'s are (resp.) $\QB$- and $L$-spread maps.
Writing $\Psi^{K/\QB}$ and $\Psi^{K/L}$ for the top and bottom compositions,
the Leray filtrations on $\underline{H_{\D}}$ for $\bar{\sX}\mapg\ms$
and $X\times\T\mapgt\T$ induce (via the two $\Psi$'s) filtrations
$\L_{K/\QB}^{\bb}\subseteq\L_{K/L}^{\bb}$ on $CH^{p}(X/K)$. (That
$\iota^{*}(\L_{\pi}^{\bb})\subseteq\L_{\pi_{\T}}^{\bb}$ follows from
Lewis's explicit description of $\L^{\bb}$ on $\underline{H_{\D}^{*}}$
on the level of representative Deligne-homology cochains.) Hence we
have maps\[
\Psi_{i}^{K/\QB}:\,\L_{(K/\QB)}^{i}CH^{p}(X/K)\to Gr_{\L_{\pi}}^{i}\underline{H_{\D}^{2p}}(\sX,\QQ(p)),\]
\[
\Psi_{i}^{K/L}:\,\L_{K/L}^{i}CH^{p}(X/K)\to Gr_{\L_{\pi_{\T}}}^{i}\underline{H_{\D}^{2p}}(X\times\eta_{\T},\QQ(p)).\]
The former extends the $\Psi_{i}^{K/\QB}$ defined {[}for the case
$\sX=X\times\gs${]} in $\S2.4$ above, but is difficult to compute;
the latter is easy to compute with $(2)$, $(3)$, $(4)$, and $\Psi_{i}^{K/L}=\iota^{*}\circ\Psi_{i}^{K/\QB}$
on $\Z\in\L_{K/\QB}^{i}$. We state what we will use:

\begin{lem}
If $\Z\in\L_{K/L}^{i}CH^{p}(X/K)$ and $\Psi_{i}^{K/L}(\Z)\neq0$,
then $\Z\nrateq0$. More precisely, one of two things is true:

(i) $\Z\notin\L_{K/\QB}^{i}$,

(ii) $\Z\in\L_{K/\QB}^{i}$ and $\Psi_{i}^{K/\QB}(\Z)\neq0$.
\end{lem}
\begin{rem}
For example $3$ it will also be helpful to note that $\Psi_{i}^{K/\QB}(\Z)$
resp. $\Psi_{i}^{K/L}(\Z)$ {}``split'' as before into $[\sZ]_{i}$
and $[AJ(\sZ)]_{i-1}$, resp. $[\sZ_{\T}]_{i}$ and $[AJ(\sZ_{\T})]_{i-1}$,
with e.g. $[\sZ]_{i}\in\hm\left(\QQ(-p),\,\underline{H^{i}}(\gs,R^{2p-i}\pi_{*}\QQ)\right)$
being sent to $[\sZ_{\T}]_{i}$ by $\iota^{*}$.
\end{rem}

\subsection{Exterior products of cycles respect the filtration}

In the more general context where $X$ may not be defined $/\QB$,
and for cycles of any codimension, it is proved in \cite{L1} that
on $CH^{*}(X/\CC)$, $\L^{i}\cdot\L^{j}\subseteq\L^{i+j}$ under the
intersection product. Moreover, push-forwards and pullbacks preserve
$\L^{\bb}$. This leads immediately to the following.

\begin{lem}
Given $\Z_{i}\in\L_{K_{i}/\QB}^{j_{i}}CH_{0}(Y_{i}/K_{i})$ for $i=1,2$,
we have $\Z:=\Z_{1}\times\Z_{2}\in\L_{K/\QB}^{j_{1}+j_{2}}CH_{0}(X/K)$
where $K:=\QB(K_{1},K_{2})$ and $X=Y_{1}\times Y_{2}$.
\end{lem}
\begin{proof}
Writing $\pi_{i}:\, X\twoheadrightarrow Y_{i}$ ($i=1,2$), $\Z=(\pi_{1}^{*}\Z_{1})\cdot(\pi_{2}^{*}\Z_{2})$
and we use the $2$ properties of $\L^{\bb}$ just mentioned.
\end{proof}
\begin{rem}
We emphasize that the $Y_{i}$ need not be defined $/\QB$, and that
$K_{1}$, $K_{2}$ need not be {}``algebraically independent'' in
the sense of Lemma $1(b)$; they could even be the same field (so
that $K=K_{1}=K_{2}$).
\end{rem}

\section{\textbf{General $\times$ special}}

In this section we present various results and examples which are
all variations on the following theme: that the product of a $0$-cycle
which is {}``general'' in some appropriate sense by a {}``special''
$0$-cycle with nontrivial Albanese image, is not rationally equivalent
to zero.

\begin{thm}
Consider $Y_{1}$ and $Y_{2}$ smooth projective varieties $/\QB$
with resp. dimensions $d_{1}$ and $d_{2}$. Let $K\subseteq\CC$
be f.g. $/\QB$ of trdeg. $j$, and $\V\in\L^{j}CH_{0}(Y_{1}/K)$
be such that its complete spread $\bar{\sV}\in Z^{d_{1}}(\ms\times Y_{1}/\QB)$
induces a nontrivial map $\o^{j}(Y_{1})\to\o^{j}(\ms)$ of holomorphic
forms. Take $\W\in CH_{0}^{\text{{hom}}}(Y_{2}/\QB)$ with $0\neq AJ_{Y_{2}}(\W)\in J^{d_{2}}(Y_{2}):=J^{d_{2}}(H^{2d_{2}-1}(Y_{2}))$.

Then $\Z:=\V\times\W\in CH_{0}(Y_{1}\times Y_{2}/K)$ is not zero
(i.e., mod $\rateq$). In particular, $\Z\in\L^{j+1}CH_{0}$ and $[\sZ]=0$;
but $[AJ(\sZ)]_{j}\neq0$ (equiv. $AJ_{X}^{j}(\Z)\neq0$), hence $\Z\notin\L^{j+2}$.
\end{thm}
\begin{proof}
The fundamental class $[\bar{\sV}]$ has $j^{\text{{th}}}$ K\"unneth
component\[
[\bar{\sV}]_{j}\in H^{j}(\ms)\otimes H^{2d_{1}-j}(Y_{1}),\]
which also lies in the r.h.s. of \small\[
\left\{ H^{j}(\ms,\CC)\otimes H^{2d_{1}-j}(Y_{1},\CC)\right\} ^{(d_{1},d_{1})}=\bigoplus_{i=0}^{j}H^{i,\, j-i}(\ms,\CC)\otimes H^{d_{1}-i,\, d_{1}+i-j}(Y_{1},\CC).\]
\normalsize We may write uniquely $[\bar{\sV}]_{j}=\sum_{i=1}^{j}[\bar{\sV}]_{(i,\, j-i)}$
and $[\bar{\sV}]_{(j,0)}=\linebreak\sum_{\ell}\a_{\ell}\otimes\nu_{\ell}\left(=\overline{[\bar{\sV}]_{(0,j)}}\right)$
where $\{\a_{\ell}\}$ is a basis of $H^{j,0}(\ms,\CC)$ and $\nu_{\ell}\in H^{d_{1}-j,\, d_{1}}(Y_{1},\CC)$.
Our hypothesis on $\bar{\sV}$ implies at least that one $\nu_{\ell}$,
say $\nu_{1}$, is nonzero.

Set $d=d_{1}+d_{2}$, $X=Y_{1}\times Y_{2}$. Since $\W$ is defined
$/\QB$, the $\QB$-spread of $\Z$ is the restriction $\sZ$ of\[
\bz=\bar{\sV}\times\W\in Z^{d}(\ms\times X/\QB)\]
to $\gs\times Y_{1}\times Y_{2}$. Already $\bz\homeq0$ with bounding
chain

\begin{equation} \d^{-1} \bz := \bar{\sV} \times \d^{-1} \W . \end{equation}${}$\\
of real dimension $2j+1$. (Here $\d^{-1}\W$ is a real $1$-chain
bounding on $\W$ which is $fixed$ for the remainder of the proof.)
Since $\V\in\L^{j}$, $[\bar{\sV}]_{i}\mapsto0\in\underline{H^{i}}(\gs)\otimes H^{2p-i}(Y_{1})$
for all $i<j$. Together with $(5)$ this implies $\left\langle \int_{\d^{-1}\bz}(\cdot)\right\rangle \mapsto0\in J^{p}(\underline{H^{i}}(\gs)\otimes H^{2p-i-1}(X))$,
hence $[AJ(\sZ)]_{i}=0$ ($\forall$ $i<j$) and $\Z\in\L^{j+1}$.

Let $[\bar{\sV}]_{j}^{\vee}\in H^{j}(\ms)\otimes H^{j}(Y_{1})$ be
any rational type $(j,j)$-class dual to%
\footnote{i.e., under {[}the restriction of{]} the Poincar\'e duality pairing
$H^{2d_{1}}(\ms\times Y_{1})\otimes H^{2j}(\ms\times Y_{1})\to\QQ$,
$[\bar{\sV}]_{j}\otimes[\bar{\sV}]_{j}^{\vee}\mapsto1$.%
} $[\bar{\sV}]_{j}$, and let $ann([\bar{\sV}]_{j})\subseteq H^{j}(\ms)\otimes H^{j}(Y_{1})$
be the subHS annihilated by $[\bar{\sV}]_{j}$. Consider the following
relative dual pairs of HS's (the $\H$'s of weight $2d-1$, the $\K$'s
of weight $2j+1$):\[
\H_{0}=H^{j}(\ms)\otimes H^{2d-j-1}(X)\supseteq H^{j}(\ms)\otimes H^{2d_{1}-j}(Y_{1})\otimes H^{2d_{2}-1}(Y_{2})=\H_{1}\,,\]
\[
\K_{0}=H^{j}(\ms)\otimes H^{j+1}(X)\supseteq H^{j}(\ms)\otimes H^{j}(Y_{1})\otimes H^{1}(Y_{2})=\K_{1}\,;\]
and \[
\H_{\sV}=\QQ[\bar{\sV}]_{j}\otimes H^{2d_{2}-1}(Y_{2})\subseteq\H_{1},\,\,\,\K_{\sV}=\QQ([\bar{\sV}]_{j}^{\vee})\otimes H^{1}(Y_{2})\subseteq\K_{1}.\]
If we take $\tilde{\Xi}=\left\langle \int_{\d^{-1}\bz}(\cdot)\right\rangle \in\left\{ F^{j+1}\K_{0}^{\CC}\right\} ^{\vee}$
and $\tilde{\xi}=[\bar{\sV}]_{j}\otimes\left\langle \int_{\d^{-1}\W}(\cdot)\right\rangle \in\left\{ F^{j+1}\K_{\sV}^{\CC}\right\} ^{\vee}$
then clearly the former annihilates $F^{j+1}\K_{2}=\linebreak F^{j+1}\left\{ ann([\bar{\sV}]_{j})_{\CC}\otimes H^{1}(Y_{2},\CC)\right\} $
completely while agreeing with $\tilde{\xi}$ on $F^{j+1}\K_{\sV}^{\CC}=\CC[\bar{\sV}]_{j}^{\vee}\otimes H^{1,0}(Y_{2},\CC)$.
Noting that $0\neq\xi=[\bar{\sV}]_{j}\otimes AJ_{Y_{2}}(\W)\in\CC[\bar{\sV}]_{j}\otimes J^{d_{2}}(Y_{2})\cong J^{d}(\H_{\sV})$
by assumption, we have from $\S2.2$ that $\iota_{\sV}(\xi)=\pr_{1}(\Xi)$.

Now let \[
\G_{0}=N^{1}H^{j}(\ms)\otimes H^{2d-j-1}(X)+H^{j}(\ms)\otimes F_{h}^{d-j}H^{2d-j-1}(X)\,\,\subseteq\,\,\H_{0},\]
so that $[AJ(\sZ)]_{j}=\b_{0}(\Xi)\subseteq J^{d}(\H_{0}/\G_{0})=J^{d}(\underline{H^{j}}(\gs)\otimes\frac{H^{2d-j-1}(X)}{F_{h}^{d-j}}).$
We must check that the image of $\G_{0}$ under the K\"unneth projection
$\widetilde{\pr_{1}}:\,\H_{0}\twoheadrightarrow\H_{1}$ is contained
in $\H_{1}\cap\G_{0}$. Set\[
\N=N^{1}H^{j}(\ms)\otimes H^{2d_{1}-j}(Y_{1})\otimes H^{2d_{2}-1}(Y_{2}),\]
\[
\F=H^{j}(\ms)\otimes F_{h}^{d-j}\left\{ H^{2d_{1}-j}(Y_{1})\otimes H^{2d_{2}-1}(Y_{2})\right\} .\]
Since $\theta:\, H^{2d-j-1}(X)\twoheadrightarrow H^{2d_{1}-j}(Y_{1})\otimes H^{2d_{2}-1}(Y_{2})$
is a morphism of HS,\[
\theta\left(F_{h}^{d-j}H^{2d-j-1}(X)\right)\,\,\,\,\subseteq\,\,\,\, F_{h}^{d-j}\left\{ H^{2d_{1}-j}(Y_{1})\otimes H^{2d_{2}-1}(Y_{2})\right\} \]
\[
\subseteq\left\{ H^{2d_{1}-j}(Y_{1})\otimes H^{2d_{2}-1}(Y_{2})\right\} \cap F_{h}^{d-j}H^{2d-j-1}(X)\]
\[
\subseteq\theta\left(F_{h}^{d-j}H^{2d-j-1}(X)\right).\]
It follows that $\widetilde{\pr_{1}}(\G_{0})\subseteq\N+\F\subseteq\H_{1}\cap\G_{0}\subseteq\widetilde{\pr_{1}}(\G_{0})$,
which gives us what we want but also that $\H_{1}\cap\G_{0}=\N+\F$.

Now by $\S2.2$ (Case 1) we are done if we can show that $\H_{\sV}\cap(\H_{1}\cap\G_{0})$
is zero in $\H_{1}$.

Take $\p:\, H^{j}(\ms,\CC)\twoheadrightarrow H^{j,0}(\ms,\CC)\oplus H^{0,j}(\ms,\CC)$
to be the projection with kernel $H^{j-1,1}(\ms,\CC)\oplus\cdots\oplus H^{1,j-1}(\ms,\CC)\supseteq N^{1}H^{j}(\ms,\CC)$.
Then the induced\[
\p_{1}:\,\H_{1}^{\CC}\twoheadrightarrow\left\{ H^{j,0}(\ms,\CC)\oplus H^{0,j}(\ms,\CC)\right\} \otimes H^{2d_{1}-j}(Y_{1},\CC)\otimes H^{2d_{2}-1}(Y_{2},\CC)\]
 kills $\N$.

Next let $\un$, $\uf$, $[\bar{\sV}]_{j}\otimes\Gamma$ be arbitrary
elements of $\N_{\CC}$, $\F_{\CC}$, and $\H_{\sV}^{\CC}$ ($\Gamma\in H^{2d_{2}-1}(Y_{2},\CC)$
arbitrary), and suppose $\un+\uf=[\bar{\sV}]_{j}\otimes\Gamma$ in
$\H_{1}^{\CC}$. In order to prove $\H_{\sV}^{\CC}\cap(\N_{\CC}+\F_{\CC})=\{0\}$,
we must show $\Gamma=0$. Apply $\p_{1}$ to both sides to get $\p_{1}(\uf)=\p_{1}([\bar{\sV}]_{j}\otimes\Gamma)$,
i.e.%
\footnote{here the bars over $\a$ and $\nu$ denote complex conjugation. %
}\[
\sum_{\ell}\a_{\ell}\otimes A_{\ell}+\sum_{\ell}\overline{\a_{\ell}}\otimes B_{\ell}=\sum_{\ell}\a_{\ell}\otimes\nu_{\ell}\otimes\Gamma+\sum_{\ell}\overline{\a_{\ell}}\otimes\overline{\nu_{\ell}}\otimes\Gamma\]
for unique classes $A_{j},\, B_{j}\in\left(F_{h}^{d-j}\left\{ H^{2d_{1}-j}(Y_{1})\otimes H^{2d_{2}-1}(Y_{2})\right\} \right)\otimes\CC$.
Hence we must have $\nu_{1}\otimes\Gamma=A_{1}$, $\overline{\nu_{1}}\otimes\Gamma=B_{1}$
in $\mathsf{{H}}_{\CC}:=H^{2d_{1}-j}(Y_{1},\CC)\otimes H^{2d_{2}-1}(Y_{2},\CC)$.
That is, $\nu_{1}\otimes\Gamma$ and $\overline{\nu_{1}}\otimes\Gamma$
must belong to {[}a subspace of{]}\[
\mathsf{{H}}_{\CC}^{d-j,d-1}\oplus\cdots\oplus\mathsf{{H}}_{\CC}^{d-1,d-j}.\]
Since $\nu_{1}$ is nonzero of pure type $(d_{1}-j,d_{1})$, $\Gamma$
belongs to \[
H^{d_{2},\, d_{2}-1}(Y_{2},\CC)\oplus\left[H^{d_{2}+1,\, d_{2}-2}(Y_{2},\CC)\oplus\cdots\oplus H^{d_{2}+j-1,\, d_{2}-j}(Y_{2},\CC)\right]\]
(bracketed terms $=0$); but since $\overline{\nu_{1}}$ is of type
$(d_{1},d_{1}-j)$, $\Gamma$ is in \[
\left[H^{d_{2}-j,\, d_{2}+j-1}(Y_{2},\CC)\oplus\cdots\oplus H^{d_{2}-2,\, d_{2}+1}(Y_{2},\CC)\right]\oplus H^{d_{2}-1,d_{2}}(Y_{2},\CC).\]
Hence $\Gamma=0$.
\end{proof}
Here is one way to construct $\V$ on $Y_{1}$ with the properties
assumed in the statement of Theorem $1$. Consider some $j$-dimensional
(possibly singular) subvariety $S_{/\QB}\subseteq Y_{1}$ with desingularization
$\ms\mapio Y_{1}$, such that the restriction ($\iota^{*}$) induces
a nontrivial map of holomorphic $j$-forms. Let $p_{0}\in\ms(\CC)$
be very general, and write $p$ for its image in $Y_{1}$. Suppose
$[\Delta_{\ms}]\in H^{j}(\ms\times\ms)$ has algebraic K\"unneth
components (also defined $/\QB$) $[\Delta_{\ms}]_{i}=[\Delta_{\ms}(i,\,2j-i)]$,
as is the case if $\ms$ is a curve, surface, abelian variety, smooth
complete intersection (in $\PP^{N}$), or arbitrary product of these.
Then $\Delta_{\ms}(j,\, j)_{*}p_{0}\in\L^{j}CH_{0}(\ms_{\CC})$, and
$\L^{\bb}$ is preserved under $\iota_{*}$ (see \cite{L1}); hence
$\V:=\iota_{*}\left[\Delta_{\ms}(j,\, j)_{*}p_{0}\right]$ lies in
$\L^{j}CH_{0}(Y_{1}^{\CC})$.

\begin{example}
If $Y_{1}\subseteq\PP^{N}$ is a (smooth) complete intersection with
$\deg(K_{Y_{1}})\geq0$, take $j=d_{1}$ and $\ms=S=Y_{1}$. The construction
gives $\V=p-o$, where $p\in Y_{1}(\CC)$ is very general and $o\in Y_{1}(\QB)$.
\end{example}
${}$

\begin{example}
Let $X=\C_{1}\times\cdots\times\C_{n}$ be a product of curves defined
$/\QB$ each of genus $\geq1$. On each $\C_{i}$ let $\W_{i}$ be
a divisor defined $/\QB$ with $AJ(\W_{i})\neq0$ in $J^{1}(\C_{i})$,
$p_{i}\in\C_{i}(\CC)$ be very general, $o_{i}\in\C_{i}(\QB)$. Assume
$\{ p_{1},\ldots,p_{n}\}$ {}``algebraically independent'' in the
sense that $p_{1}\times\cdots\times p_{n}\in X(\CC)$ is very general
(not defined over a field of trdeg $<n$). Then for each $j\geq1$,
take\[
S=\C_{1}\times\cdots\times\C_{j}\times\{ o_{j+2}\}\times\cdots\times\{ o_{n}\}\subseteq\C_{1}\times\cdots\times\C_{j}\times\C_{j+2}\times\cdots\times\C_{n}=Y_{1}\,,\]
$p_{0}=p_{1}\times\cdots\times p_{j}(\times o_{j+2}\times\cdots\times o_{n})\in S(\CC)$.
One obtains from the construction \[
\V=(p_{1}-o_{1})\times\cdots\times(p_{j}-o_{j})\times o_{j+2}\times\cdots\times o_{n}\in\L^{j}CH_{0}(Y_{1}^{\CC})\]
 which has $cl_{X}^{j}(\V)\neq0$; hence it follows from the Theorem
that \[
\Z=(p_{1}-o_{1})\times\cdots\times(p_{j}-o_{j})\times\W_{j+1}\times o_{j+2}\times\cdots\times o_{n}\in\L^{j+1}CH_{0}(X_{\CC})\]
has $AJ_{X}^{j}(\Z)\neq0$. Thus $\V,\,\Z\,\nrateq0$.

In \cite{K1}, this example is tied to Bloch's results (\cite{B1})
for $0$-cycles on Abelian varieties (which are dominated by these
products of curves). One also gets applications to Calabi-Yau 3-folds. 
\end{example}
Here is a slight generalization of Theorem $1$ that gives some initial
evidence that our invariants are well-behaved under products.

\begin{cor}
Let $Y_{1},\, Y_{2},\,\W$ be as in Theorem $1$ and take $\tilde{\V}\in\L^{j}CH_{0}(Y_{1}^{\CC})$
with $cl_{Y_{1}}^{j}(\tilde{\V})\neq0$. Then if the GHC holds, $AJ_{X}^{j}(\tilde{\V}\times\W)\neq0$
hence $\tilde{\V}\times\W\nrateq0$.
\end{cor}
\begin{proof}
We have $\tilde{K}$ finitely generated $/\QB$ and $\tilde{\ms}$
(defined $/\QB$), such that $\tilde{\V}/\tilde{K}$ and $\QB(\tilde{\ms})\cong\tilde{K}$.
By definition of $cl_{Y_{1}}^{j}$, $0\neq[\tv]_{j}\in\linebreak\hm\left(\QQ(-d_{1}),\,\underline{H^{j}}(\eta_{\tilde{\ms}})\otimes H^{2d_{1}-j}(Y_{1})\right)$.
By Lemma $3(b),(c)$ (and HC), $\exists$ a $j$-dimensional section
$\ms_{/\QB}\mapiota\tilde{\ms}$ such that (for $\bv=\iota^{*}\btv$)
$[\sV]_{j}=\iota^{*}[\tv]_{j}\neq0$. By GHC$(1,j,\ms)$ and Lemma
$2$, $\bv^{*}:\,\o^{j}(Y_{1})\to\o^{j}(\ms)$ is nontrivial hence
$[AJ(\sV\times\W)]_{j}\neq0$ by the Theorem. Now obviously $\tilde{\V}\times\W\in\L^{j+1}$
and $[\btv\times\W]=0$, hence $[AJ(\tv\times\W)]_{j}$ is defined
and maps to $[AJ(\sV\times\W)]_{j}$ under $\iota^{*}$.
\end{proof}
In the Theorem and Corollary $Y_{1},\, Y_{2},\,\W$ are defined $/\QB$,
and we would like to have a more flexible statement. One idea is to
allow $Y_{1}$ (like $\V$) to be defined $/K$, in such a way that
the generalized $[\sV]_{j}$ of Remark $2$ is still $\neq0$. However
this turns out to be too optimistic as we now show.

\begin{example}
\textbf{(Bielliptic cycle)} 

Let $E_{\l}:=\overline{\{ y^{2}=x(x-1)(x-\l)\}}\mapm\PP_{[x]}^{1}$;
these form a family $\E\mapg\M\subseteq\PP_{[\l]}^{1}$. Writing $o_{\l}:=(o,o)$,
pick $\l_{2}\in\QB$ and $q_{\l_{2}}\in E_{\l_{2}}(\QB)$ s.t. $AJ(q_{\l_{2}}-o_{\l_{2}})\in J^{1}(E_{\l_{2}})$
is nonzero. Put $\e:=\mathsf{{X}}_{\l_{2}}(q_{\l_{2}})\in\PP^{1}(\QB)$
and choose $q_{\l}\in\mathsf{{X}}_{\l}^{-1}(\e)$ (with $\e$ fixed)
continuously in $\l\in\CC$, which requires lifting to a double cover
$\tilde{\E}\to\tilde{\M}$ of the family.

Define $\Z_{\l}=q_{\l}-o_{\l}\in CH_{0}(E_{\l})$, take $\l_{1}\in\CC\setminus\QB$
and set $K:=\QB(\l_{1})$. The cycle $\Z_{\l_{1}}$ on $E_{\l_{1}}$
is defined over {[}an algebraic extension of{]} $K$; its spread $\sZ_{1}$
on $\tilde{\E}$ yields a normal function $\nu\in\Gamma(\tilde{\M},\,\J_{\tilde{\E}/\tilde{\M}}^{1})$
defined by $\nu(\l):=AJ_{E_{\l}}(\Z_{\l})$. One easily sees this
is nontrivial (since it is $2$-torsion at $\l=\e$ but not elsewhere);
arguing by monodromy, its infinitesimal invariant $\delta\nu\in\Gamma\left(\tilde{\M},\,\frac{\o_{\tilde{\M}}^{1}\otimes\H_{\tilde{\E}/\tilde{\M}}^{1}}{\nabla(\O_{\tilde{\M}}\otimes\F^{1}\H_{\tilde{\E}/\tilde{\M}}^{1}}\cong\o_{\tilde{\M}}^{1}\times\H_{\tilde{\E}/\tilde{\M}}^{1,0}\right)$
must be nonzero.

According to Remark $2$, the generalized $\Psi_{1}^{K/\QB}(\Z_{\l_{1}})$
invariant maps to a generalized $[\sZ_{1}]_{1}\in\hm\left(\QQ(-1),\underline{H^{1}}(\eta_{\tilde{\M}},R^{1}\pi_{*}\QQ)\right)$
which itself maps (injectively) to the infinitesimal invariant, hence
$[\sZ_{1}]_{1}\neq0$. So in Theorem $1$'s notation, $\Z_{\l_{1}}$
plays the role of $\V$ (in a generalized sense), and $\Z_{\l_{2}}$
(defined $/\QB$, with nontrivial $AJ$ image) that of $\W$.

But $\Z_{\l_{1}}\times\Z_{\l_{2}}=(q_{\l_{1}}-o_{\l_{1}})\times(q_{\l_{2}}-o_{\l_{2}})\rateq0$
(modulo $2$-torsion). This is because $(q_{\l_{1}},q_{\l_{2}})$
is a (very general) point on the image of the bielliptic curve $E_{\l_{1}}\times_{\PP_{[x]}^{1}}E_{\l_{2}}=:B$
in $E_{1}\times E_{2}$. Writing $\sigma$ for $(-\text{{id}})$ on
each $E_{i}$, $\sigma\times\sigma$ induces the hyperelliptic involution
on $B$; hence an explicit function on $B$ gives \[
-(q_{\l_{1}},q_{\l_{2}})-(\sigma(q_{\l_{1}}),\sigma(q_{\l_{2}}))+2(o_{\l_{1}},o_{\l_{2}})\rateq0.\]
One now easily shows $-(\sigma\times\sigma)(\Z_{\l_{1}}\times\Z_{\l_{2}})\rateq\Z_{\l_{1}}\times\Z_{\l_{2}}\rateq(\sigma\times\sigma)(\Z_{\l_{1}}\times\Z_{\l_{2}})\,;$
the assertion follows.
\end{example}
Here, then, is the appropriate generalization.

\begin{cor}
Let $K\supseteq L$ be an extension (of trdeg $\geq j$) of subfields
of $\CC$ f.g. $/\QB$. Let $Y_{1},\, Y_{2},\,\W$ be defined $/L$
but otherwise as in Theorem $1$. Referring to $\S2.5$, take $\V\in\L_{K/\QB}^{j}CH_{0}(Y_{1}/K)$
s.t. $\bar{\sV_{\T}}\in Z^{d_{1}}\left((\T\times Y_{1})/L\right)$
induces $\o^{j}(Y_{1})\to\o^{j}(\T)$ nontrivial.

Then $\Z:=\V\times\W\nrateq0$; more precisely, $\Psi_{j+1}^{K/\QB}(\Z)\neq0$.
\end{cor}
\begin{proof}
By the proof of Theorem $1$, we have $0\neq[AJ(\sZ_{\T})]_{j}^{tr}\in\linebreak J^{d}\left(\underline{H^{j}}(\eta_{\T})\otimes\frac{H^{2d-j-1}(X)}{F_{h}^{d-j}}\right).$
Hence, $\Psi_{j+1}^{K/L}(\Z)\neq0$. Since $\V\in\L_{K/\QB}^{j}$
and $\W\in\L_{K/\QB}^{1}$, $\Z\in\L_{K/\QB}^{j+1}$ by Lemma $5$;
hence by Lemma $4$, $\Psi_{j+1}^{K/\QB}(\Z)\neq0$.
\end{proof}
\begin{rem}
We expect no better, in the sense that there are situations where
$[\sZ_{\T}]_{j+1}=0$, $[AJ(\sZ_{\T})]_{j}\neq0$ for the $partial$
spread but $[\sZ]_{j+1}\neq0$ for the $\QB$-spread. (See \cite[sec. 7.1]{K1}.)
\end{rem}

\section{\textbf{General $\times$ general}}

As far as taking products of cycles that $both$ spread is concerned,
the easiest case is where each has a nontrivial higher cycle-class.
Let $Y_{1},\, Y_{2}$ be defined $/\QB$. If $\V_{i}\in\L^{j_{i}}CH_{0}(Y_{i}/K_{i})$
(for $i=1,2$) have $cl_{Y_{i}}^{j_{i}}(\V_{i})\neq0$, then with
{}``reasonable'' assumptions we can expect $cl_{Y_{1}\times Y_{2}}^{j_{1}+j_{2}}(\V_{1}\times\V_{2})\neq0$
(clearly $\V_{1}\times\V_{2}\in\L^{j_{1}+j_{2}}$ by Lemma $5$).
Namely, we need a hypothesis that guarantees the spread of $\V_{1}\times\V_{2}$
to be just $\bv_{1}\times\bv_{2}$ on $(\ms_{1}\times Y_{1})\times(\ms_{2}\times Y_{2})$;
algebraic {}``independence'' of $K_{1}$ and $K_{2}$ in the sense
of Lemma $1(b)$ is sufficient. Applying the GHC and Lemma $2$ in
each factor, each $\bv_{i}$ induces a nonzero map $\o^{j_{i}}(Y_{i})\to\o^{j_{i}}(\ms_{i})$;
hence $\bv_{1}\times\bv_{2}$ does the same $\o^{j_{1}+j_{2}}(Y_{1}\times Y_{2})\to\o^{j_{1}+j_{2}}(\ms_{1}\times\ms_{2})$
and so $[\sV_{1}\times\sV_{2}]_{j_{1}+j_{2}}\neq0$. If we make no
assumption on $K_{1}$ and $K_{2}$ then there are problems:

\begin{example}
Let $Y_{1}=Y_{2}=E\,/\QB$ be an elliptic curve and $\V_{1}=\V_{2}=p-o$
where $p\in E(\CC)$ is very general and $o$ is a $\QB$-point. Then
$[\sV_{1}]_{1}\neq0$, $[\sV_{2}]_{1}\neq0$ but $\V_{1}\times\V_{2}=(p,p)-(o,p)-(p,o)+(o,o)$
is the diagonal cycle; this is $\rateq0$ (mod $2$-torsion).
\end{example}
\begin{rem}
Should one want to generalize to $Y_{1},\, Y_{2}$ not defined $/\QB$,
the above assumptions --- nontriviality of $cl_{Y_{i}}^{j_{i}}(\V_{i})$
($i=1,2$) and {}``independence'' of the fields of definition of
$\V_{1}$ and $\V_{2}$ (namely, $K_{1}$ and $K_{2}$) over $\QB$
--- are insufficient to guarantee $cl_{Y_{1}\times Y_{2}}^{j_{1}+j_{2}}(\V_{1}\times\V_{2})\neq0$.
For a counterexample (in case $j_{1}=j_{2}=1$) one can simply take
$\l_{1}$ and $\l_{2}$ both general (and algebraically independent
over $\QB$) in Example $3$. However, if one takes $K_{1}$ and $K_{2}$
independent over the common field of defintion of $Y_{1}$ and $Y_{2}$
(say, $L$), then an analogous result obviously holds for the higher
$cl$-type invariants arising from the partial $L$-spreads of $\V_{1}$,
$\V_{2}$, and $\V_{1}\times\V_{2}$.
\end{rem}
Now suppose (with $Y_{i}$ def'd. $/\QB$) $\V_{i}\in\L^{j_{i}}CH_{0}(Y_{i}/K_{i})$
but $[\sV_{i}]_{j_{i}}=0$, $[AJ(\sV_{i})]_{j_{i}-1}\neq0$ ($i=1,2$);
i.e. each cycle has nontrivial higher $AJ$-class. The situation looks
more grim here for nontriviality of the exterior product, even if
we assume $K_{1}$ and $K_{2}$ independent. It will never be the
case that $[AJ(\sV_{1}\times\sV_{2})]_{j_{1}+j_{2}-2}\neq0$, because
$\V_{1}\times\V_{2}\in\L^{j_{1}+j_{2}}$ by Lemma $5$. In fact, if
we started with $\trdeg(K_{i})=j_{i}-1$ ($i=1,2$) then all the higher
cycle- and $AJ$-classes of $\V_{1}\times\V_{2}$ are zero; if an
extension of the Bloch-Beilinson conjecture holds (see \cite{K2}
or \cite{L1}) then this $\implies$ $\rateq0$.

The interesting problem is the asymmetric one: $\V_{1}$ with nontrivial
higher $cl$, $\V_{2}$ with higher $AJ$-class $\neq0$. Referring
to Example $2$, if we take $\V_{1}=(p_{1}-o_{1})\times\cdots\times(p_{m}-o_{m})$
(where $m<n$) and $\V_{2}=(p_{m+1}-o_{m+1})\times\cdots\times(p_{n-1}-o_{n-1})\times\W_{n}$
on $Y_{1}=\C_{1}\times\cdots\times\C_{m}$ and $Y_{2}=\C_{m+1}\times\cdots\times\C_{n}$
(resp.), then of course the product cycle is nontrivial. More generally,
one expects any $\Z$ from Theorem $1$ to work as $\V_{2}$.

Here is the stongest general result we could prove; note $\V$, $\W$
replace $\V_{1}$, $\V_{2}$. Naturally we would have preferred to
assume only (say) $AJ_{Y_{2}}^{j_{2}}(\tilde{\W})^{tr}\neq0$; see
Remark $6$ for a conditional improvement along these lines.

\begin{thm}
Let $Y_{1},\, Y_{2}$ $/\QB$ be smooth projective w./dimensions $d_{1},\, d_{2}$;
$\tilde{\V}\in\L^{j_{1}}CH_{0}(Y_{1}/\CC)$ with $cl_{Y_{1}}^{j_{1}}(\tilde{\V})\neq0$;
$\tilde{\W}\in\L^{j_{2}+1}CH_{0}(Y_{2}/\CC)$ with $\overline{AJ_{Y_{2}}^{j_{2}}(\tilde{\W})^{tr}}\neq0$
and $cl_{Y_{2}}^{j_{2}+1}(\tilde{\W})=\cdots=cl_{Y_{2}}^{d_{2}}(\tilde{\W})=0$.
Assume $\tilde{\V},\,\tilde{\W}$ have resp. models over fields $K_{1},\, K_{2}\subseteq\CC$
f.g. $/\QB$ with trdegs. $t_{1},\, t_{2}$, such that $\QB(K_{1},K_{2})$
has trdeg. $t_{1}+t_{2}$. Assume the GHC.

Then $\tilde{\V}\times\tilde{\W}\in\L^{j_{1}+j_{2}+1}CH_{0}(Y_{1}\times Y_{2}/\CC)$
has $AJ_{Y_{1}\times Y_{2}}^{j_{1}+j_{2}}(\tilde{\V}\times\tilde{\W})^{(tr)}\neq0$.
\end{thm}
\begin{proof}
Let $\bar{\tilde{\sW}}$ be a choice of {}``complete'' spread. Since
$cl_{Y_{2}}^{i}(\tilde{\W})=0$ for $0\leq i\leq d_{2}$, the image
$[\tilde{\sW}]_{i}\in\hm\left(\QQ(-d_{2}),\,\underline{H^{i}}(\eta_{\tilde{\ms}_{2}})\otimes H^{2d_{2}-i}(Y_{2})\right)$
of $[\bar{\tilde{\sW}}]_{i}$ is zero for all $i$ (automatic for
$i>d_{2}$ since $Gr_{\L}^{i}CH^{d_{2}}=0$ by \cite{L1}). Thus $[\bar{\tilde{\sW}}]_{i}\in N^{1}H^{i}(\tilde{\ms}_{2})\otimes H^{2d_{2}-i}(Y_{2})$
and by Deligne \cite[Cor. 8.2.8]{D} there exist irreducible codim.-$1$
$\QB$-subvarieties $S_{\a}$ on $\tilde{\ms}_{2}$ such that $[\bar{\tilde{\sW}}]_{i}\in\text{{Gy}}\left\{ \oplus_{\a}H^{i-2}(\tilde{S}_{\a})\otimes H^{2d_{2}-i}(Y_{2})\right\} .$
Hence $[\bar{\tilde{\sW}}]$ is a sum of Gysin images of classes in
$\hm\left(\QQ(-d_{2}+1),\, H^{2d_{2}-2}(\tilde{S}_{\a}\times Y_{2})\right)$.
By the HC, these are given by cycles; thus one may modify $\bar{\tilde{\sW}}$
(without affecting $\tilde{\sW}$) so that $[\bar{\tilde{\sW}}]=0$.

By Lemma $1(b)$, the (complete) spread of $\tilde{\V}\times\tilde{\W}$
is just the product of spreads, $\bar{\tilde{\sV}}\times\bar{\tilde{\sW}}$.
Now we specialize in both factors as in Lemma $3$ and Remark $1$,
obtaining (with HC) $\V$ and $\W$ exactly as in the hypotheses of
the Proposition below. (We also have to use GHC($1,j_{1},\ms_{1}$)
to get from $[\sV]_{j_{1}}\neq0$ to the map of holomorphic forms.)
According to the Proposition, $[AJ(\sV\times\sW)]_{j_{1}+j_{2}}^{tr}\neq0$;
and $[AJ(\tilde{\sV}\times\tilde{\sW})]_{j_{1}+j_{2}}^{tr}$ maps
to this under the specialization.%
\footnote{That this map is well-defined simply follows from well-definedness
of $H^{j_{1}+j_{2}}(\tilde{\ms}_{1}\times\tilde{\ms}_{2})/N^{1}\,\,\to\,\, H^{j_{1}+j_{2}}(\ms_{1}\times\ms_{2})/N^{1}$
(argue as in proof of Lemma $3(b)$).%
} 
\end{proof}
So the proof has been reduced to this statement, which is what we
could prove $without$ assuming GHC.

\begin{prop*}
Let $\V\in\L^{j_{1}}CH_{0}(Y_{1}/k_{1})$, $\W\in\L^{j_{2}+1}CH_{0}(Y_{2}/k_{2})$
for $k_{1},\, k_{2}\subseteq\CC$ of resp. $\trdeg_{/\QB}$ $j_{1}$
and $j_{2}$. Assume that $\bar{\sV}$ induces a nontrivial map of
holomorphic $j_{1}$-forms $\o^{j_{1}}(Y_{1})\to\o^{j_{1}}(\ms_{1})$;
that $\overline{[AJ(\sW)]_{j_{2}}^{tr}}\neq0$; and that $\W$ has
a complete spread $\bar{\sW}$ with $[\bar{\sW}]=0$. Then $[AJ(\sV\times\sW)]_{j_{1}+j_{2}}^{(tr)}\neq0$.
\end{prop*}
\begin{proof}
\emph{(of Proposition)} Clearly the first paragraph of the proof of
Theorem $1$ applies (replacing $\ms$, $j$ by $\ms_{1}$, $j_{1}$).
Define $[\bar{\sV}]_{j}$, $[\bar{\sV}]_{j}^{\vee}$, $ann\left([\bar{\sV}]_{j}\right)$
as before. Set $d=d_{1}+d_{2}$, $j=j_{1}+j_{2}$, $\ms=\ms_{1}\times\ms_{2}$,
$X=Y_{1}\times Y_{2}$, $\bar{\sZ}=\bar{\sV}\times\bar{\sW}$, and
$\d^{-1}\bar{\sZ}=\bar{\sV}\times\d^{-1}\bar{\sW}$.

We begin with the HS (omitting the obvious dual $\K$'s)\[
\H_{0}\,=\, H^{j}(\ms)\otimes H^{2d-j-1}(X)\,\,\,\,\,\,\supseteq\]
\[
H^{j_{1}}(\ms_{1})\otimes H^{j_{2}}(\ms_{2})\otimes H^{2d_{1}-j_{1}}(Y_{1})\otimes H^{2d_{2}-j_{2}-1}(Y_{2})\,=\,\H_{1}\]
 and \[
\G_{0}=F_{h}^{1}H^{j}(\ms)\otimes H^{2d-j-1}(X)+H^{j}(\ms)\otimes F_{h}^{d-j}H^{2d-j-1}(X).\]
Reasoning as in the proof of Thm. $1$, under $\H_{0}\twoheadrightarrow\H_{1}$,
$\G_{0}$ projects to\[
\G_{0}\cap\H_{1}\,=\, F_{h}^{1}\left\{ H^{j_{1}}(\ms_{1})\otimes H^{j_{2}}(\ms_{2})\right\} \otimes H^{2d_{1}-j_{1}}(Y_{1})\otimes H^{2d_{2}-j_{2}-1}(Y_{2})\]
\[
\mspace{100mu}+\, H^{j_{1}}(\ms_{1})\otimes H^{j_{2}}(\ms_{2})\otimes F_{h}^{d-j}\left\{ H^{2d_{1}-j_{1}}(Y_{1})\otimes H^{2d_{2}-j_{2}-1}(Y_{2})\right\} \]
\[
=:\,\N\,+\,\F.\mspace{290mu}\]

Now define $\tilde{\Xi}=\left\langle \int_{\d^{-1}\bar{\sZ}}(\cdot)\right\rangle \in\left\{ F^{j+1}\K_{0}^{\CC}\right\} ^{\vee}$,
so that $[AJ(\sZ)]_{j}^{tr}$ projects to%
\footnote{equals if GHC($1,j,\ms$) holds%
} $\b_{0}(\Xi)\in J^{d}(\H_{0}/\G_{0})$; this is what we must show
nonzero. Write $\sigma_{23}:\,\H_{1}\mapf H^{j_{1}}(\ms_{1})\otimes H^{2d_{1}-j_{1}}(Y_{1})\otimes H^{j_{2}}(\ms_{2})\otimes H^{2d_{2}-j_{2}-1}(Y_{2})$
for the map exchanging the $2^{\text{{nd}}}$ and $3^{\text{{rd}}}$
$\otimes$-factors, and define $\H_{\sV}$ by \[
\sigma_{23}(\H_{\sV})=\QQ[\bar{\sV}]_{j_{1}}\otimes H^{j_{2}}(\ms_{2})\otimes H^{2d_{2}-j_{2}-1}(Y_{2}).\]
Set $\tilde{\xi}=[\bar{\sV}]_{j_{1}}\otimes\left\langle \int_{\d^{-1}\bar{\sW}}(\cdot)\right\rangle \in\left\{ F^{j+1}\K_{\sV}^{\CC}\right\} ^{\vee}$;
then $\tilde{\Xi}$ annihilates $F^{j+1}\K_{2}^{\CC}$ and agrees
with $\tilde{\xi}$ on $F^{j+1}\K_{\sV}^{\CC}$, so that $\iota_{\sV}(\xi)=\pr_{1}(\Xi)$.
Now define $\G_{1}$ by \[
\sigma_{23}(\G_{1})=H^{j_{1}}(\ms_{1})\otimes H^{2d_{1}-j_{1}}(Y_{1})\otimes SF_{h}^{(1,\, d_{2}-j_{2})}\left\{ H^{j_{2}}(\ms_{2})\otimes H^{2d_{2}-j_{2}-1}(Y_{2})\right\} .\]
By assumption (nontriviality of $\overline{[AJ(\sW)]_{j_{2}}^{tr}}$),
$\pi_{\sV}(\xi)$ is nonzero in$\linebreak$$J\left(\frac{\H_{\sV}}{\G_{1}\cap\H_{\sV}}\right)\cong\CC[\bar{\sV}]_{j}\otimes J^{d_{2}}\left(\frac{H^{j_{2}}(\ms_{2})\otimes H^{2d_{2}-j_{2}-1}(Y_{2})}{SF_{h}^{(1,\, d_{2}-j_{2})}\{\text{{num}}\}}\right).$
According to $\S2.2$ (Case $2$) we are done modulo showing\[
(\N+\F)\cap\H_{\sV}\subseteq\G_{1}.\]

Projecting along the Hodge decomposition in each factor,\small \[
\p_{1}'{}:\,\H_{1}^{\CC}\twoheadrightarrow H^{j_{1},0}(\ms_{1},\CC)\otimes H^{j_{2},0}(\ms_{2},\CC)\otimes H^{d_{1}-j_{1},d_{1}}(Y_{1},\CC)\otimes H^{d_{2}-j_{2}-1,d_{2}}(Y_{2},\CC)\]
 \normalsize kills $\N$ and $\F$. Write $\{\a_{i}\}$ for a basis
of $H^{j_{1},0}(\ms_{1},\CC)$, and $\{\Gamma_{\ell}\}_{\ell=1}^{M}$
for a basis of $H^{2d_{2}-j_{2}-1}(Y_{2},\CC)$ s.t. $\{\Gamma_{\ell}\}_{\ell=1}^{r}\subseteq H^{d_{2}-j_{2}-1,d_{2}}(Y_{2},\CC)$
and $\{\Gamma_{\ell}\}_{\ell=r+1}^{M}\subseteq F^{d_{2}-j_{2}}H^{2d_{2}-j_{2}-1}(Y_{2},\CC)$.
Let\[
\un+\uf\,=\,\sum_{\ell=1}^{M}\sigma_{23}^{-1}\left([\bar{\sV}]_{j_{1}}\otimes\g_{\ell}\otimes\Gamma_{\ell}\right)\]
be an arbitrary element of $\left(\N_{\CC}+\F_{\CC}\right)\cap\H_{\sV}^{\CC}$,
where each $\g_{\ell}\in H^{j_{2}}(\ms_{2},\CC)$ has a Hodge decomposition
$\sum_{p+q=j_{2}}\g_{\ell}^{(p,q)}$. Applying $\p_{1}'$ gives\[
0\,=\,\sum_{i}\sum_{\ell=1}^{r}\a_{i}\otimes\g_{\ell}^{(j_{2},0)}\otimes\nu_{i}\otimes\Gamma_{\ell}\,;\]
since $\nu_{1}\neq0$ this implies \[
0\,=\,\sum_{\ell=1}^{r}\a_{1}\otimes\g_{\ell}^{(j_{2},0)}\otimes\nu_{1}\otimes\Gamma_{\ell}\,.\]
Hence for $\ell=1,\ldots,r$ we have $\g_{\ell}^{(j_{2},0)}=0$, i.e.
$\g_{\ell}\in\overline{F^{1}H^{j_{2}}(\ms_{2},\CC)}$.%
\footnote{here (and similarly below) the bar denotes complex conjugation%
} It follows that

\begin{equation} \begin{matrix} \sum_{\ell =1}^M \g_{\ell} \otimes \Gamma_{\ell} \mspace{25mu} \in \mspace{25mu} \overline{F^1 H^{j_2} (\ms_2 , \CC )} \otimes H^{2d_2 -j_2 -1} (Y_2 ,\CC ) \\ \mspace{220mu} + \mspace{10mu} H^{j_2}(\ms_2, \CC ) \otimes F^{d_2 - j_2} H^{2d_2 -j_2 -1} (Y_2, \CC ) . \end{matrix} \end{equation}${}$\\
An argument symmetric in Hodge types (starting from the definition
of $\p_{1}'$, replace all $(p,q)$'s by $(q,p)$'s) shows

\begin{equation} \begin{matrix} \sum_{\ell =1}^M \g_{\ell} \otimes \Gamma_{\ell} \mspace{25mu} \in \mspace{25mu} F^1 H^{j_2} (\ms_2 , \CC ) \otimes H^{2d_2 -j_2 -1} (Y_2 ,\CC ) \\ \mspace{220mu} + \mspace{10mu} H^{j_2}(\ms_2, \CC ) \otimes \overline{F^{d_2 - j_2} H^{2d_2 -j_2 -1} (Y_2, \CC )} . \end{matrix} \end{equation}${}$\\
Hence $\sum\g_{\ell}\otimes\Gamma_{\ell}$ lives in the intersection\[
F^{1}H^{j_{2}}(\ms_{2},\CC)\otimes F^{d_{2}-j_{2}}H^{2d_{2}-j_{2}-1}(Y_{2},\CC)\mspace{100mu}\]
\[
\mspace{100mu}+\,\overline{F^{1}H^{j_{2}}(\ms_{2},\CC)}\otimes\overline{F^{d_{2}-j_{2}}H^{2d_{2}-j_{2}-1}(Y_{2},\CC)}\]
\[
=\, SF^{(1,\, d_{2}-j_{2})}\left\{ H^{j_{2}}(\ms_{2},\CC)\otimes H^{2d_{2}-j_{2}-1}(Y_{2},\CC)\right\} \]
of $(6)$ and $(7)$; and so $(\N+\F)\cap\H_{\sV}$ lies in {[}the
$\CC$-vector space{]} $\sigma_{23}^{-1}\left(\QQ[\bar{\sV}]_{j_{1}}\otimes SF^{(1,d_{2}-j_{2})}\right)$.
But $(\N+\F)\cap\H_{\sV}$ is a HS, hence it actually lies in the
largest subHS (of, say, $\H_{\sV}$) contained in $\linebreak$$\sigma_{23}^{-1}\left(\QQ[\bar{\sV}]_{j_{1}}\otimes SF^{(1,d_{2}-j_{2})}\right)$,
which (using purity of $\QQ[\bar{\sV}]_{j_{1}}$) is {[}the $\QQ$-vector
space{]} $\sigma_{23}^{-1}\left(\QQ[\bar{\sV}]_{j_{1}}\otimes SF_{h}^{(1,d_{2}-j_{2})}\right)$.
Obviously this gives containment in $\G_{1}$ and completes the proof.
\end{proof}
\begin{rem}
(i) Theorem $1$ is the case $j_{2}=0$.

(ii) It is easy to show (by sharpening slightly the argument in the
proof of Thm. $1$) that the cycles $\Z$ with nontrivial $[AJ(\sZ)]_{j}^{tr}$
produced by Thm. $1$ actually have nontrivial $\overline{[AJ(\sZ)]_{j}^{tr}}$,
hence would make a suitable choice of $\W$ for the above Proposition.
(If this were not the case, one would expect a stronger Proposition
to be true!)

(iii) Assume the GHC. Then in the statement of the above Proposition,
we may relax the requirement on $\W$ to $[AJ(\sW)]_{j_{2}}^{tr}\neq0$
provided\[
F_{h}^{1}\left\{ H^{j_{1}}(\ms_{1})\otimes H^{j_{2}}(\ms_{2})\right\} =F_{h}^{1}H^{j_{1}}(\ms_{1})\otimes H^{j_{2}}(\ms_{2})+H^{j_{1}}(\ms_{1})\otimes F_{h}^{1}H^{j_{2}}(\ms_{2})\]
 and \[
F^{d-j}\left\{ H^{2d_{1}-j_{1}}(Y_{1})\otimes H^{2d_{2}-j_{2}-1}(Y_{2})\right\} \,\,\,=\]
\[
F_{h}^{d_{1}-j_{1}+1}H^{2d_{1}-j_{1}}(Y_{1})\otimes H^{2d_{2}-j_{2}-1}(Y_{2})+H^{2d_{1}-j_{1}}(Y_{1})\otimes F_{h}^{d_{2}-j_{2}}H^{2d_{2}-j_{2}-1}(Y_{2}).\]
(These are decidedly $not$ satisfied if e.g. $\ms_{1}=\ms_{2}$.)
\end{rem}
\begin{proof}
\emph{(of Rem. 6(iii))} In this case $\N+\F$ is annihilated by the
projection from $\H_{1}$ to\[
\frac{H^{j_{1}}(\ms_{1})}{F_{h}^{1}}\otimes\frac{H^{j_{2}}(\ms_{2})}{F_{h}^{1}}\otimes\frac{H^{2d_{1}-j_{1}}(Y_{1})}{F_{h}^{d_{1}-j_{1}+1}}\otimes\frac{H^{2d_{2}-j_{2}-1}(Y_{2})}{F_{h}^{d_{2}-j_{2}}},\]
and it follows that $\sum\g_{\ell}\otimes\Gamma_{\ell}$ lies in $F_{h}^{1}H^{j_{2}}(\ms_{2})\otimes H^{2d_{2}-j_{2}-1}(Y_{2})+H^{j_{2}}(\ms_{2})\otimes F_{h}^{d_{2}-j_{2}}H^{2d_{2}-j_{2}-1}(Y_{2})$.
Writing $\sigma_{23}(\G_{1}'{})$ for $H^{j_{1}}(\ms_{1})\otimes H^{2d_{1}-j_{1}}(Y_{1})$
tensor this, $\G_{1}'$ replaces $\G_{1}$ in the above argument and
$\pi_{\sV}(\xi)$ need only be assumed nontrivial in \[
J^{d}\left(\frac{\H_{\sV}}{\G_{1}'{}\cap\H_{\sV}}\right)\,\,\cong\,\,\CC[\bar{\sV}]_{j_{1}}\otimes J^{d_{2}}\left(\underline{H^{j_{2}}}(\eta_{\ms_{2}})\otimes\frac{H^{2d_{2}-j_{2}-1}(Y_{2})}{F_{h}^{d_{2}-j_{2}}\{\text{{num}}\}}\right)\]
(assuming GHC($1,j_{2},\ms_{2}$) for the $\cong$).
\end{proof}

\address{\noun{Max Planck Institut f\"ur Mathematik, Vivatsgasse 7,}}

\address{\noun{53111 Bonn, Germany}}

\email{kerr@mpim-bonn.mpg.de}
\end{document}